\documentclass[12pt]{article}
\usepackage{amsmath,enumerate, amsfonts, amssymb, amsthm, graphicx}
\usepackage{epic,eepic}
\def\UseSection{
        \numberwithin{equation}{section}
    \theoremstyle{plain}
        \newtheorem{theorem}    {Theorem}[section]
        \DefineTheorems 
}

\def\DefineTheorems{
    
    \newtheorem{lemma}      [theorem] {Lemma}
    
    \newtheorem{prop}       [theorem] {Proposition}
    
    \newtheorem{cor}        [theorem] {Corollary}

    \theoremstyle{definition}
    \newtheorem{defn}       [theorem] {Definition}

    \theoremstyle{definition}

}

\newcommand{\bt}   {\begin{theorem}}
\newcommand{\et}   {\end  {theorem}}
\newcommand{\bl}   {\begin{lemma}}
\newcommand{\el}   {\end  {lemma}}
\newcommand{\bp}   {\begin{prop}}
\newcommand{\ep}   {\end  {prop}}
\newcommand{\bc}   {\begin{cor}}
\newcommand{\ec}   {\end  {cor}}
\newcommand{\bd}   {\begin{defn}}
\newcommand{\ed}   {\end  {defn}}

\newcommand{\ba}   {\begin{array}}
\newcommand{\ea}   {\end  {array}}
\newcommand{\be}   {\begin{enumerate}}
\newcommand{\ee}   {\end  {enumerate}}
\newcommand{\bi}   {\begin{itemize}}
\newcommand{\ei}   {\end  {itemize}}

\def\eq#1\en{\begin{equation}#1\end{equation}}
\def\eqsplit#1\ensplit{
    \begin{equation}\begin{split}#1\end{split}\end{equation}
    }
\def\eqalign#1\enalign{
    \begin{align}#1\end{align}
    }
\def\eqmul#1\enmul{
    \begin{multline}#1\end{multline}
    }
\newcommand{\eqarrstar} {\begin{eqnarray*}}
\newcommand{\enarrstar} {\end{eqnarray*}}
\newcommand{\eqarray}   {\begin{eqnarray}}
\newcommand{\enarray}   {\end{eqnarray}}

\newcommand{\lbeq}[1]  {\label{e:#1}}
\newcommand{\refeq}[1] {\eqref{e:#1}}    

\makeatletter
\newcommand{\labelcounter}[2]{{%
    \stepcounter{#1}
    \protected@write\@auxout{}%
    {\string\newlabel{#2}{{\csname the#1\endcsname}{\thepage}}}%
    {\ref{#2}}
    }}
\makeatother
\newcommand{\sss}   { \scriptscriptstyle }

\newcommand{\Ebold} {{\mathbb E}}

\newcommand{\Pbold} {{\mathbb P}}

\newcommand{\Tbold} {{\mathbb T}}

\newcommand{\Zbold} {{\mathbb Z}}

\newcommand{\Ccal}   {\mathcal{C}}

\newcommand{\Zd}    {{ {\Zbold}^d }}

\newcommand{\spose}[1] {{\hbox to 0pt{#1\hss}} }
\newcommand{\ltapprox} {\mathrel{\spose{\lower 3pt\hbox{$\mathchar"218$}}
 \raise 2.0pt\hbox{$\mathchar"13C$}}}
\newcommand{\gtapprox} {\mathrel{\spose{\lower 3pt\hbox{$\mathchar"218$}}
 \raise 2.0pt\hbox{$\mathchar"13E$}}}

\UseSection  
\newcommand{\vep}{\varepsilon}
\newcommand{\Z}{\mathbb{Z}}
\newcommand{\sZ}{{\sss \mathbb{Z}}}
\newcommand{\sT}{{\sss \mathbb{T}}}
\newcommand{\R}{\mathbb{R}}
\newcommand{\torus}{\mathbb T_{r,d}}
\newcommand{\T}{\mathbb T}
\newcommand{\Br}{B_{r,d}}
\newcommand{\gr}{{\mathbb G}}

\newcommand{\egr}\Exp
\newcommand{\cn}{\Omega}
\newcommand{\bb}{\underline{b}}
\newcommand{\tb}{\overline{b}}
\newcommand{\shift}{\!\!\!\!}
\newcommand{\Cmax}{{\cal C}_{\rm max}}
\newcommand{\CZ}{{\cal C}_{\sZ}}
\newcommand{\Ctorus}{{\cal C}_{\sT}}
\newcommand{\Cg}{C_{\sss\rm\chi}}
\newcommand{\cxi}{c_{\rm{\sss \xi}}}
\newcommand{\Cxi}{C_{\rm{\sss \xi}}}
\newcommand{\Ctc}{C_{\rm{\tilde\chi}}}

\newcommand{\Exp}{{{\mathbb E}_p}}

\newcommand{\conn}{\longleftrightarrow}
\newcommand{\nn}{\nonumber}

\newcommand{\simr}{\stackrel{r}{\sim}}
\newcommand{\Tconn}{\stackrel{\sT}{\conn}}
\newcommand{\Zconn}{\stackrel{\sZ}{\conn}}

\begingroup\makeatletter\ifx\SetFigFont\undefined%
\gdef\SetFigFont#1#2#3#4#5{%
  \reset@font\fontsize{#1}{#2pt}%
  \fontfamily{#3}\fontseries{#4}\fontshape{#5}%
  \selectfont}%
\fi\endgroup%

\title  {
        Random graph asymptotics on high-dimensional tori
        }

\author{
Markus Heydenreich$~^*$
\and
Remco van der Hofstad \footnote{Department of Mathematics and
Computer Science, Eindhoven University of Technology, P.O.~Box~513,
5600~MB Eindhoven, The~Netherlands.
E-mail: {\tt m.o.heydenreich@tue.nl, r.w.v.d.hofstad@tue.nl}}
}

\begin{document}
\maketitle

\begin{abstract}
We investigate the scaling of the largest critical
percolation cluster on a large $d$-di\-men\-sional torus,
for nearest-neighbor percolation in sufficiently high dimensions,
or when $d>6$ for sufficiently spread-out percolation.
We use a relatively simple coupling
argument to show that this largest critical cluster is,
with high probability, bounded above by a large constant times $V^{2/3}$ and
below by a small constant times $V^{2/3}(\log{V})^{-4/3}$,
where $V$ is the volume of the torus. We also give
a simple criterion in terms of the subcritical percolation two-point function on
$\Z^d$ under which the lower bound can be improved to small constant
times $V^{2/3}$, i.e.\ we prove
{\it random graph asymptotics} for the largest
critical cluster on the high-dimensional torus. This establishes
a conjecture by \cite{Aize97}, apart from logarithmic corrections.
We discuss implications of these results
on the dependence on boundary conditions for high-dimensional
percolation.

Our method is crucially based on the results in \cite{BCHSS04a, BCHSS04b},
where the $V^{2/3}$ scaling was proved {\it subject} to the assumption that
a suitably defined critical window contains the percolation threshold
on $\Z^d$. We also strongly rely on mean-field results for percolation on
$\Z^d$ proved in \cite{Hara90, Hara05, HHS03, HS90a}.
\end{abstract}

\section{Introduction}
\subsection{The model}
We consider Bernoulli bond percolation on the graph $\gr$, where $\gr$ is either the hypercubic lattice $\Zd$, or the finite torus $\torus=\{-\lfloor r/2\rfloor,\dots,\lceil r/2\rceil-1\}^d$.

For $\gr=\Zd$, we consider two sets of bonds. In the {\em nearest-neighbor model}, two vertices $x$ and $y$ are linked by a bond whenever $|x-y|=1$, whereas in the {\em spread-out model}, they are linked whenever $0<\|x-y\|\le L$. Here, and throughout the paper, we write $\|\cdot\|$ for the supremum norm, and $|\cdot|$ for the Euclidean norm on $\Zd$. The integer parameter $L$ is typically chosen large. We let each bond independently be occupied with probability $p$, or vacant otherwise. The resulting product measure is denoted by $\Pbold_{\sZ,p}$, and the corresponding expectation $\Ebold_{\sZ,p}$.
We write $\{0\conn x\}$ for the event that there exists a path of occupied bonds from the origin $0$ to the lattice site $x$, and define
    \eq
    \lbeq{tau-def}
    \tau_{\sZ,p}(x):=\Pbold_{\sZ,p}(0\conn x)
    \en
to be the \emph{two-point} function.
We further write $\CZ(x):=\{y\in\Zd\mid x\conn y\}$ for the \emph{cluster} or connected component of $x$, $|\CZ(x)|$ for the number of vertices in $\CZ(x)$ and $\chi_\sZ(p):=\sum_{x\in\Zd}\tau_{\sZ,p}(x)=\Ebold_{\sZ,p}|\CZ(0)|$ for the expected cluster size.
The degree of the graph, which we denote by $\cn$, is thus $\cn=2d$ in the nearest-neighbor case and $\cn=(2L+1)^d-1$ in the spread-out case.

It is well-known that bond percolation on $\Zd$ in dimension $d\ge2$ obeys a phase transition, i.e.\ there exists a critical threshold $p_c(\Zd)\in(0,1)$ such that $p_c(\Zd)=\inf\{p\colon\Pbold_{\sZ,p}(|\CZ(0)|=\infty)>0\}$. Furthermore, by the results in \cite{AB87,Men86}, $p_c(\Zd)$ can be expressed as $p_c(\Zd)=\sup\{p\colon\chi_\sZ(p)<\infty\}$.

For $\gr=\torus$, we also consider two related settings:
\begin{enumerate}
\item The nearest-neighbor torus:
an edge joins vertices that differ by $1$ (modulo $r$)
in exactly one component.
For $d$ fixed and $r$ large, this is a periodic approximation to $\Z^d$.
Here $\Omega =2d$ for $r \geq 3$.
We study the limit in which $r\rightarrow \infty$ with $d>6$ fixed, but large.

\item The spread-out torus:
an edge joins vertices $x=(x_1,\ldots,x_d)$ and $y=(y_1,\ldots,y_d)$
if $0<\max_{i=1,\ldots,d}|x_i-y_i|_r \leq L$ (with $|\cdot|_r$ the metric
on $\Zbold_r$).  We study the limit $r \to \infty$,
with $d > 6$ fixed and $L$ large (depending on $d$)
and fixed.  This gives a periodic approximation to
range-$L$ percolation on $\Z^d$. Here $\Omega = (2L+1)^d-1$
provided that $r\geq 2L+1$, which we will always assume.
\end{enumerate}
We consider bond percolation on these tori with bond occupation probability $p$ and write $\Pbold_{\sT,p}$ and
$\Ebold_{\sT,p}$ for the product measure and corresponding expectation, respectively.
We use the notation $\tau_{\sT,p}(\cdot)$, $\chi_\sT(p)$ and $\Ctorus(\cdot)$ analogously to the corresponding $\Zd$-quantities.

In this paper, we will investigate the size of the maximal cluster on $\torus$, i.e.,
    \eq
    |\Cmax|:=\max_{x\in\torus}|\Ctorus(x)|,
    \en
at the critical percolation threshold $p_c(\Zd)$.
An alternative definition for the critical percolation threshold on the torus, denoted by $p_c(\torus)$, was given in \cite[(1.7)]{BCHSS04a} as the solution to
    \eq\lbeq{defPcT}
    \chi_{\sT}(p_c(\torus))=\lambda V^{1/3},
    \en
where $\lambda$ is a sufficiently small constant, and $V=|\torus|=r^d$ denotes the volume of the torus.
The definition of $p_c(\torus)$ in \refeq{defPcT} is an \emph{internal} definition only, due to the fact that \cite{BCHSS04b} deals with rather general tori, for which an external definition (such as $p_c(\Zd)$) does not always exist.
On the other hand, the internal definition in \refeq{defPcT} assumes \emph{a priori} mean-field behaviour, and is therefore unsuitable outside this setting. On the high-dimensional torus $\torus$, we therefore have \emph{two} sensible critical values, the externally defined $p_c(\Zd)$, and the internally defined $p_c(\torus)$ in \refeq{defPcT}.
One of the goals of this paper is to investigate how close these two critical values are.

The most prominent example of percolation on a finite graph is the \emph{random graph}, which is obtained by applying percolation to the \emph{complete} graph. This has been first studied by Erd\H{o}s and R\'{e}nyi in 1960 \cite{ER60}. They showed that, when $p$ is scaled as $(1+\vep)V^{-1}$, there is a phase transition at $\vep=0$. For $\vep<0$, the size of the largest cluster is proportional to $\log V$, whereas for $\vep>0$, it is proportional to $V$.
For $\vep=0$, the size of the largest cluster divided by $V^{2/3}$ weakly converges to some (non-trivial) limiting random variable, while the expected cluster size is, as in \refeq{defPcT}, proportional to $V^{1/3}$.
This follows from results by Aldous \cite{Ald97}.
See also \cite{Bol85} for results up to 1984, and \cite{JKL93,JLR00,LPW94} for references to subsequent work.
We will refer to the $V^{2/3}$-scaling as \emph{random graph asymptotics}.

In this paper, we study the size of the largest cluster on the torus for $p=p_c(\Zd)$. It has been shown by Borgs, Chayes, van der Hofstad, Slade and Spencer \cite{BCHSS04a,BCHSS04b} that, if
    \eq\lbeq{scalingWindow}
    p_c(\Zd)=p_c(\torus)+O(V^{-1/3}),
    \en
then, with probability at least $1-O(\omega^{-1})$, $|\Cmax|$ is in between $\omega^{-1}V^{2/3}$ and $\omega V^{2/3}$ as $V\to\infty$, for $\omega\ge1$ sufficiently large.
Here, we write $f(x)=O(g(x))$ for functions $f,g\ge0$ and $x$ converging to some limit, if there exists a constant $C>0$ such that $f(x)\le Cg(x)$ in the limit, and $f(x)=o(g(x))$ if $g(x)\neq O(f(x))$.
Furthermore, we write  $f=\Theta(g)$ if $f=O(g)$ and $g=O(f)$.

Aizenman \cite{Aize97} conjectured that this random graph asymptotics holds for the maximal critical cluster in dimension $d>6$, as we explain in more detail below.
Also in \cite{BCHSS04b} it was conjectured that \refeq{scalingWindow} holds.
By means of a coupling argument, we prove that a slightly weaker statement than \refeq{scalingWindow} (with a logarithmic correction in the lower bound, see \refeq{pcZdbd} below) indeed holds for $d$ sufficiently large in the nearest-neighbor model, or $d>6$ and $L$ sufficiently large in the spread-out model.
Furthermore, we give a criterion which we believe to hold, and which implies \refeq{scalingWindow} without logarithmic corrections.

Note that all our results assume that $d$ is large in the nearest-neighbor model or $d>6$ and $L$ large in the spread-out model. That is, we require the torus to be in some sense high-dimensional. We do believe that the results hold for all $d>6$ and $L\ge1$, however, the proof relies on various lace expansion results, which require that the degree $\cn$ is large. On the other hand, we do not expect these asymptotics to be true for $d\le6$.

Aizenman \cite{Aize97} studied a similar question, but now for percolation
on a box of width $r$ under {\it bulk} boundary conditions, where clusters are defined to be the intersection of the box $\{-\lfloor r/2\rfloor,\dots,\lceil r/2\rceil-1\}^d$ with clusters in the infinite lattice (and thus clusters need not to be connected within the box).
Aizenman assumed that the probability, at criticality, that $x$ is connected to the origin
is bounded above and below by constants times $\|x\|^{-(d-2)}$. This assumption was
established in \cite{HHS03}
for the spread-out model for $d>6$ and sufficiently large, but finite $L\ge1$,
and in \cite{Hara05} for the nearest-neighbor
model above 19 dimensions.
Aizenman showed that, under this condition on the two-point function, the size of the largest connected component under bulk boundary conditions is, with high probability, bounded from above by a constant times $r^4\log r$, and bounded from below by $\vep_r\,r^4$ for any sequence $\vep_r\to0$ as $r\to\infty$.
Furthermore, he conjectures that the $r^4$-scaling for the size of the largest cluster holds for dimension $d>6$ also under \emph{free} boundary conditions (where no connections outside the box are allowed), but changes to $V^{2/3}= r^{2d/3}\gg r^4$ under \emph{periodic} boundary conditions.
This indicates the importance of boundary conditions at criticality in high dimensions.
We will further elaborate on the role of boundary conditions in Section \ref{sec-boundaryConditions}.

\subsection{Results}\label{subsect-results}
Our first result gives asymptotic bounds on the size of the largest cluster.
\begin{theorem}
\label{thm-1}
Fix $d>6$ and $L$ sufficiently large in the spread-out case, or $d$ sufficiently large for nearest-neighbor
percolation.
Then there exist constants $b_1, b_2, C>0$, such that for all $\omega_1\ge C$ and $\omega_2\ge1$,
    \eq
    \lbeq{Cmaxbd}
    \Pbold_{\sT, p_c(\Z^d)}\Big(\omega_1^{-1}V^{2/3}(\log{V})^{-4/3}
    \leq |\Cmax|\leq \omega_2 V^{2/3}\Big)
    \geq 1-\frac{b_1}{\omega_1^{3/2}(\log V)^2}-\frac{b_2}{\omega_2}
    \qquad\text{as $r\to\infty$.}
    \en
The constant $b_1$ can be chosen as $288\cdot120^{3}$, and $b_2$ equal to $b_6$ in \cite[Theorem 1.3]{BCHSS04a}.
\end{theorem}

To prove Theorem \ref{thm-1}, we will use a coupling argument relating $\chi_\sT(p)$ and $\chi_\sZ(p)$ to show that there exists a constant $\Lambda\ge0$ such that, when $r\rightarrow \infty$,
    \eq
    \lbeq{pcZdbd}
    p_c(\torus)-\frac{\Lambda}{\cn} V^{-1/3}(\log{V})^{2/3}
    \leq p_c(\Z^d)\leq p_c(\torus)+\frac{\Lambda}{\cn} V^{-1/3}.
    \en
Relying on results in \cite{BCHSS04a}, \refeq{pcZdbd} implies \refeq{Cmaxbd}.
Inequality \refeq{Cmaxbd} implies that $|\Cmax|V^{-2/3}$ is a tight random variable, but it does not rule out that $|\Cmax|V^{-2/3}\to0$ as $V\to\infty$.

Our method is crucially based on the results in \cite{BCHSS04a, BCHSS04b},
but we also rely on mean-field results for percolation on
$\Z^d$ by Hara \cite{Hara90, Hara05}, Hara and Slade \cite{HS90a}, and Hara, van der Hofstad and Slade \cite{HHS03}.
Each of these papers relies on the \emph{lace expansion}, but the lace expansion will not be used in this paper.
The lecture notes by Slade \cite{Sla04} and Hara and Slade \cite{HS94} provide a general introduction to the lace expansion and its role in proving mean-field critical behavior for percolation and related models.

Note that in \cite{BCHSS04a,BCHSS04b}, $p_c(\torus)$ was defined as in \refeq{defPcT}.
The results in \cite{BCHSS04a,BCHSS04b}, however, do not establish rigorously that the exponent $1/3$ in \refeq{defPcT} is the only correct choice.
Indeed, \cite{BCHSS04a,BCHSS04b} suggest that a smaller exponent would also do, since the supercritical results proved there are not sufficiently sharp.
Theorem \ref{thm-1} shows that, at least in terms of the power of $V$, the scaling of $|\Cmax|$ at $p_c(\Zd)$ and at $p_c(\torus)$ is identical, thus establishing that on $\Zd$ the choice \refeq{defPcT} \emph{is} appropriate.

Unfortunately, the lower bound in Theorem \ref{thm-1} does not quite meet the
upper bound. Under a condition on the percolation two-point function, we can prove
the matching lower bound.
To state this result, we introduce the quantity
    \eq
    \lbeq{tildechidef}
    \widetilde \chi_{\sZ}(p,r):=\sup_{y}\sum_{z\simr y,\, \|z\|
    \geq \frac r2} \tau_{\sZ,p}(z),
    \en
where we write that $x\simr y$ when $x~({\rm mod}~r)=y~({\rm mod}~r)$.
We will call $x$ and $y$ \emph{r-equivalent} when $x\simr y$.

\begin{theorem}
\label{thm-2}
Under the assumptions in Theorem \ref{thm-1},
suppose that there exists a $K>0$ such that, for
$p=p_c(\Z^d)-{K}{\cn}^{-1}V^{-1/3}$ and some $C_{\sss K}>0$, the bound
    \eq
    \lbeq{sharpassumpthm}
    \widetilde \chi_{\sZ}(p,r)\leq C_{\sss K} V^{-2/3}
    \en
holds.
Then, for all $\omega\ge1$ and $b$ equal to $b_6$ in \cite[Theorem 1.3]{BCHSS04a},
    \eq
    \lbeq{Cmaxubthm}
    \Pbold_{\sT, p_c(\Z^d)}\left(\frac1\omega V^{2/3}\le|\Cmax|\le\omega V^{2/3}\right)
    \geq 1-\frac{b}{\omega}.
    \en
\end{theorem}

\noindent
Inequality \refeq{Cmaxubthm} implies that $|\Cmax|V^{-2/3}$ is tight, and that each possible weak limit along any subsequence is non-zero.
The result in \refeq{Cmaxubthm} combined with the results in \cite{BCHSS04a} would indicate that the scaling of $|\Cmax|$ at $p_c(\torus)$ and at $p_c(\Zd)$ agree, thus showing that there is no significant difference between the internally and externally defined critical values.

Analogously to Theorem \ref{thm-1}, we will show that when \refeq{sharpassumpthm} holds, there exists a constant $\Lambda\ge0$ such that, when
$r\rightarrow \infty$,
    \eq
    \lbeq{pcZdubd}
    p_c(\torus)-\frac{\Lambda}{\cn} V^{-1/3}\leq p_c(\Z^d)\leq p_c(\torus)+\frac{\Lambda}{\cn} V^{-1/3},
    \en
and deduce \refeq{Cmaxubthm} using \cite{BCHSS04a}.

We strongly believe that \refeq{sharpassumpthm} holds.
Indeed, \refeq{sharpassumpthm} follows when the two-point function is sufficiently
smooth. For example, the condition
    \eq\lbeq{smoothTau}
    \max_{\|z\|\ge \frac r2,\,x\in\torus}\frac{\tau_{\sZ,p}(z)}{\tau_{\sZ,p}(z+x)}\le C
    \en
for some positive constant $C$,
where we consider $\torus=\{-\lfloor r/2\rfloor,\dots,\lceil r/2\rceil-1\}^d$ as a subset of $\Zd$,
implies that
    \eq
    \tau_{\sZ,p}(z)\le\frac CV\sum_{x\in\torus}\tau_{\sZ,p}(z+x)
    \qquad\text{for $\|z\|\ge \frac r2$.}
    \en
Note that, for every $y\in\Zd$, we have that
$\sum_{z\simr y}\sum_{x\in\torus}f(z+x)=\sum_{x\in\Zd}f(x)$
for all functions $f\colon \Zd\mapsto \R$.
Consequently,
    \eq
    \lbeq{tildechirestr}
    \widetilde \chi_{\sZ}(p,r)
    =\sup_{y}\sum_{z\simr y,\,\|z\|\ge \frac r2} \tau_{\sZ,p}(z)
    \le \frac CV \sup_{y}\!\!\sum_{z\simr y,\,\|z\|\ge \frac r2} \,\sum_{x\in\torus} \tau_{\sZ,p}(z+x)
    \le \frac CV \sum_{x\in\Zd} \tau_{\sZ,p}(x)
    =\frac CV\,\chi_{\sZ}(p).
    \en
Thus, for $p=p_c(\Z^d)-{K}{\cn}^{-1}V^{-1/3}$, we obtain by the fact that $\gamma=1$
(see \cite{AN84} and \cite{HS90a}, or Theorem \ref{prop-gamma} below) that $\widetilde \chi_{\sZ}(p,r)$ is bounded from above by a constant multiple of $V^{-2/3}$.

\subsection{Related results}
In this section we discuss the relation between Theorems \ref{thm-1} and \ref{thm-2} and the literature.

Hara and Slade \cite{HS00a,HS00b} study the geometry of large critical clusters on the rescaled lattice. Under the conditions of Theorem \ref{thm-1}, they show that critical clusters with size of order $n$ on the lattice rescaled by $n^{-1/4}$ converge to integrated super-Brownian excursion as $n\to\infty$.
Together with the results in \cite{HS90a}, and using \cite{AB87,AN84}, these papers prove that various critical exponents for percolation exist and take on mean-field values.

Borgs, Chayes, Kesten and Spencer \cite{BCKS99,BCKS01}
consider the largest cluster in a finite box of width $r$ under \emph{free} boundary conditions, i.e., clusters are connected only within the box.
They show that, for $p=p_c(\Zd)$, the largest critical cluster scales like $V^{\delta/(1+\delta)}$, where the critical exponent $\delta$ is defined by
    \eq
    \Pbold_{p_c(\Zd)}\big(|\Ccal_\sZ(0)|\ge n\big)\approx n^{-1/\delta}\quad\text{as $n\to\infty$},
    \en
under some conditions related to the so-called scaling and hyperscaling postulates.
The hyperscaling postulates are proven in dimension $d=2$, and are widely believed to hold up to the upper critical dimension $6$.
For the mean-field value $\delta=2$, proved in \cite{HS00a,HS00b}, we would obtain the $V^{2/3}$ asymptotics.
However, in \cite{BCKS99,BCKS01}, it was assumed that crossing probabilities of a cube of dimensions $(r,3r,\dots,3r)$ remain uniformly bounded away from $1$ as $r\to\infty$.
In high dimensions, Aizenman \cite{Aize97} proves that any cube $\{0,\dots,r\}^d$ has crossings with high probability, so that the results in \cite{BCKS99,BCKS01} do not apply. Also, in high dimensions, the hyperscaling relations are not valid.
More specifically, one hyperscaling relation is that
    \eq
    2-\eta=d\,\frac{\delta-1}{\delta+1}.
    \en
Under the conditions of Theorem \ref{thm-1}, due to \cite{BA91,Hara05,HHS03,HS90a}, we have that $\eta=0$, $\delta=2$, so that this hyperscaling relation fails for $d>6$.

Theorems \ref{thm-1} and \ref{thm-2} study the scaling of the largest critical percolation cluster on the high-dimensional torus. These results indicate that the scaling limit of the largest critical cluster should be described by $|\Cmax|V^{-2/3}$. We conjecture that, at $p=p_c(\Zd)$, the random variables $|\Cmax|V^{-2/3}$ converge as $r\to\infty$ to some (non-trivial) limiting distribution. It would be of interest to investigate whether, if the rescaled largest cluster converges, the limit law is identical to the limit of $|\Cmax|n^{-2/3}$ for the largest cluster of the random graph on $n$ vertices, as identified by Aldous \cite{Ald97}.
The convergence of $|\Cmax|V^{-2/3}$ would describe part of the {\em incipient infinite cluster} (IIC) for percolation on the torus, as described by Aizenman \cite{Aize97}.
Aizenman's IIC is closely related to the scaling limit of percolation on large cubes, see \cite[Section~5]{Aize97}. Mind also the warning at the bottom of \cite[p.~553]{Aize97}.

Another approach to the incipient infinite cluster is described by Kesten \cite{Kes86a}. Indeed, Kesten investigates the local configuration close to the origin in $\Z^2$, conditioned on the critical cluster of the origin to be infinite. Since, at criticality, the cluster of the origin is infinite with probability $0$, an appropriate limit needs to be taken. Kesten offers two alternatives:
\begin{enumerate}
\item[(i)] To condition the origin to be connected to infinity at $p>p_c(\Z^2)$ and take the limit $p\searrow p_c(\Z^2)$.
\item[(ii)] To condition the critical cluster of the origin to be connected to the boundary of the box $\{-n,\dots,0,\dots,n\}^2$ and take the limit $n\to\infty$.
\end{enumerate}
Kesten proves that both limits exist and are equal. This limit is Kesten's incipient infinite cluster.
Kesten was motivated to describe this IIC in order to study random walk on large critical clusters \cite{Kes86b}, for which physicists have performed simulations showing subdiffusive behavior.

J\'arai \cite{Jar03b,Jar03a} extended these results, and proved that several other natural conditioning and limiting schemes give the same limit. In one of these constructions, J\'arai takes a uniform point in the largest critical cluster on a box $\{0,\dots,r-1\}^d$, shifts it to the origin and takes the limit $r\to\infty$.
In \cite{HJ04}, the proof of existence of the IIC was extended to high-dimensional percolation, under the assumptions of Theorems \ref{thm-1} and \ref{thm-2}. The proof in \cite{HJ04} follows the proof in \cite{HHS02}, where the IIC was constructed for spread-out oriented percolation above 4 spatial dimensions.

We conjecture that, as $r\to\infty$, the law of local configurations around a uniform point in $\Cmax$ at criticality converges to the IIC as constructed in \cite{HJ04}. This result would give a natural link between the scaling limit of critical percolation on a large box in \cite{Aize97} and Kesten's notion of the IIC in \cite{Kes86a}.

We recall that, following the conjecture in \cite{Aize97}, the size of the largest connected component $|\Cmax|$ on the cube $\{0,\dots,r-1\}^d$ under \emph{free} bondary conditions scales like $r^4$ (as under bulk boundary conditions), in contrast to the $V^{2/3}$-scaling under periodic boundary conditions.
Such qualitatively different behavior between free and periodic boundary conditions has also been observed when studying loop-erased random walks and uniform spanning trees on a finite box in high dimensions, as we will explain now.

Choose two uniform points $x$ and $y$ from the $d$-dimensional box of side length $r$, with $d>4$.
We are interested in the graph distance between these two points on a uniform spanning tree.
Pemantle \cite{Pema91} showed that this graph distance has the same distribution as the length of a loop-erased random walk starting in $x$ and stopped when reaching $y$.
Loop-erased random walk above 4 dimensions converges to Brownian motion (cf.\ \cite[Section 7.7]{Lawl91}), that is, it scales diffusively.
This suggests that, under \emph{free} boundary conditions, the graph distance between $x$ and $y$ scales like $r^2$.
On the other hand, the combined results of Benjamini and Kozma \cite{BK05} and Peres and Revelle \cite{PR05} show that the distance between $x$ and $y$ on a uniformly chosen spanning tree on the \emph{torus} $\torus$ is of the order $V^{1/2}=r^{d/2}>r^2$ for $d>4$.
Schweinsberg \cite{Schw06} identifies the logarithmic correction for the scaling on the 4-dimensional torus.

\subsection{Organization}
This paper is organized as follows.
In Section \ref{sec-coupling}, we state a coupling result that is crucial for all our subsequent bounds.
In Section \ref{sec-previous}, we collect the main results from previous work that are used in our arguments.
In Section \ref{sec-upperBound}, we prove the upper bound in Theorem \ref{thm-1}.
In Section \ref{sec-lowerBound}, we prove the complementary lower bounds in Theorems \ref{thm-1} and \ref{thm-2}.
Finally, in Section \ref{sec-boundaryConditions}, we discuss how the growth of the maximal cluster depends on the precise boundary conditions.

\section{A coupling result for clusters on the torus and $\Z^d$}
\label{sec-coupling}
In this section, we prove that the cluster size for percolation
on the torus is stochastically smaller than the one on $\Z^d$ by
a coupling argument. We fix $p$, and omit the subscript
$p$ from the notation. We use subscripts $\Z$ and
$\T$ to denote objects on $\Z^d$ and $\torus$, respectively.

The goal of this section is to give a coupling of the $\torus$-cluster and the $\Zd$-cluster of the origin. This will be achieved by constructing
these two clusters
simultaneously from a percolation configuration on $\Zd$, as we explain in more detail now.

The basic idea is that, on any graph, it is well-known that the law of a cluster
$\Ccal(0)$ can be described by subsequently exploring the bonds one can reach
from 0.
We will first describe this exploration of a cluster in some detail, before giving the coupling, which is described by a more
elaborate way of exploring the percolation clusters on the torus and on $\Zd$ simultaneously
from a percolation configuration on $\Z^d$.
The exploration process is
defined in terms of {\it colors} of the bonds.
Initially, all bonds are uncolored, which means that they have not yet
been explored. During the exploration process we will color the bonds
black if they are found to be occupied, and white if they are found to
be vacant. Furthermore, we distinguish between active and inactive vertices.
Initially, only the origin $0$ is active, and all other vertices are inactive.

We now explore the bonds in the graph according to the following scheme.
We order the vertices in an arbitrary way. Let $v$ be the smallest active vertex.
Now we explore (and color) all uncolored bonds that have an endpoint in $v$, i.e.,
we make the bond black with probability $p$ and white with
probability $1-p$, independently of all other bonds.
In case we have assigned the black color, we set the vertex
at the other end of the bond active (unless it was already active and none of its neighboring edges are now uncolored, in which case we make it inactive).
In particular, the active vertices are those vertices that are part of a black bond, as well as an uncolored bond.
Finally, after all bonds starting at $v$ have been explored, we set $v$ inactive. We repeat doing so until there are no more active vertices.
In the latter case, the exploration process is completed, i.e., there are no more black bonds that share a common endpoint with an uncolored bond.
The cluster $\Ccal(0)$ is equal to the set of vertices that are part of the black bonds.
When the graph is finite, then this procedure always stops. When
the graph is infinite, then the exploration process continues
forever precisely when $|\Ccal(0)|=\infty$. This completes the exploration
of a single cluster on a general graph.

The exploration of a single cluster will be extended to explore the cluster on the torus $\Ctorus(0)$ and the cluster on the infinite lattice $\CZ(0)$ simultaneously from a percolation configuration on $\Zd$. For $\CZ(0)$, the result of the exploration will be identical to the exploration of a single cluster described above. The related cluster $\Ctorus(0)$ is a subset of all vertices that are $r$-equivalent to vertices part of a black bond.
The main result in this section is Proposition \ref{prop-coupling}, whose proof gives the details of this simultaneous construction of the two clusters.

To state the result, we need some notation.
For $x,y\in \torus$, we write $x\Tconn y$  when $x$
is connected to $y$ in the percolation configuration on the torus,
while, for $x,y\in \Z^d$, we write $x\Zconn y$ when $x$ is connected
to $y$ in the percolation configuration on $\Z^d$.
Also, we call two distinct bonds $\{x_1, y_1\}$ and
$\{x_2, y_2\}$ {\it $r$-equivalent} if there exists an element
$z\in \Z^d$ such that $\{x_1,y_1\}=\{x_2+r z,y_2+r z\}$.
We sometimes abbreviate $r$-equivalent to equivalent.
For a directed bond $b=(x,y)$, we write $\bb=x$ and $\tb=y$, and
for two bonds $b_1=(x_1,y_1)$ and $b_2=(x_2,y_2)$, we write $b_1\simr b_2$ when
$(x_1,y_1)=(x_2+rz,y_2+rz)$ for some $z\in\Zd$.
For $A$ and $B$ increasing events, we denote by $A\circ B$ the event that $A$ and $B$ occur on disjoint sets of bonds, see \cite{BK85}.

\bp[The coupling]
\label{prop-coupling}
Consider nearest-neighbor percolation for $r\ge3$ or spread-out percolation for $r\ge 2L+1$, in any dimension.
There exists a probability law $\Pbold_{\sZ,\sT}$ on the joint space of $\Zd$- and $\torus$-percolation such that, for all events $E$,
    \eq
    \Pbold_{\sZ,\sT}(\Ctorus(0)\in E)=
    \Pbold_{\sT}(\Ctorus(0)\in E), \qquad \Pbold_{\sZ,\sT}(\CZ(0)\in E)=
    \Pbold_{\sZ}(\CZ(0)\in E),
    \en
and $\Pbold_{\sZ,\sT}$-almost surely, for all $x\in \torus$,
    \eq
    \lbeq{subsetCs}
    \{0\Tconn x\}\subseteq  \bigcup_{y\in \Z^d\colon y\simr x}\{0\Zconn y\}.
    \en
In particular, $|\Ctorus(0)|\le|\CZ(0)|$.
Moreover, for $x\simr y$, and $\Pbold_{\sZ,\sT}$-almost surely,
    \eqalign\lbeq{lemBK}
    &
    \{0\Zconn y\}\cap \{0\Tconn x\}^c\\
    &\qquad \subseteq \shift \bigcup_{b_1\neq b_2\colon b_1\simr b_2}\bigcup_{z\in \Z^d}
    (0\Zconn z)\circ (z\Zconn \bb_1)\circ (z\Zconn \bb_2)\circ
    (b_2 \text{ is }\Z-\text{occ.})\circ (\tb_2\Zconn y).\nn
    \enalign

\ep
\vskip0.5cm

\noindent
Equation \refeq{subsetCs} will be used to conclude that the expected cluster size on $\torus$ is bounded from above by the one on $\Zd$. In order to prove our main results, we use \refeq{lemBK} to prove a related lower bound on the expected cluster size on $\torus$ in terms of the one on $\Zd$. See Sections \ref{sec-upperBound} and \ref{sec-lowerBound} for details.
The inequality on the cluster sizes of $\torus$- and $\Zd$-percolation also follows from the coupling used by Benjamini and Schramm \cite[Theorem 1]{BenjaSchra96}. However, \refeq{lemBK} does not follow immediately from their work.

\vspace{.5cm}
\noindent
{\it Proof of Proposition \ref{prop-coupling}.}
The exploration of a single cluster, as described above, will be generalized to construct $\Ctorus(0)$ and $\CZ(0)$ simultaneously
from a percolation configuration on $\Zd$.
The difference between
percolation on the torus and on $\Z^d$ can be summarised
by saying that, on the torus, $r$-equivalent bonds have the
{\it same} occupation status, while on $\Z^d$, equivalent and distinct
bonds have an {\it independent} occupation status.
For the exploration of the torus $\Ctorus(0)$, we have to make sure that we
explore equivalent bonds at most \emph{once}. We therefore
introduce a third color, gray, indicating that the bond
itself has not been explored yet, but one of its equivalent
bonds has. Therefore, at each step of the exploration process, we have 4 different types of bonds on $\Z^d$:

\begin{itemize}

\item uncolored bonds, which have not been explored yet;

\item black bonds, which have been explored and found to be
occupied;

\item white bonds, which have been explored and found to be
vacant;

\item gray bonds, of which an equivalent bond has been explored.
\end{itemize}

As in the exploration of a single cluster, we number the vertices of $\Zd$ in an arbitrary way, and start with all bonds uncolored and only the origin active.
Then we repeat choosing the smallest active bond, and explore all uncolored bonds containing it. However, after exploring a bond (and coloring
it black or white), we color all bonds that are $r$-equivalent
to it gray.
Again, this exploration is completed when there are no more
active vertices. This is equivalent to the fact that there are no
more black (and therefore occupied) bonds sharing a common endpoint with an
uncolored bond.

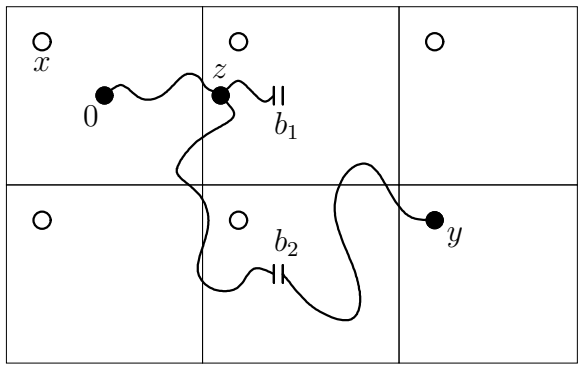
\begin{figure}[h]
\centering
{   \setlength{\unitlength}{0.03pt}
{\renewcommand{\dashlinestretch}{30}
\begin{picture}(7224,4539)(0,-10)
\put(5412,1812){\blacken\ellipse{224}{224}}
\put(5412,1812){\ellipse{224}{224}}
\put(1249,3387){\blacken\ellipse{224}{224}}
\put(1249,3387){\ellipse{224}{224}}
\put(2712,3387){\blacken\ellipse{224}{224}}
\put(2712,3387){\ellipse{224}{224}}
\thicklines
\put(462,4062){\ellipse{224}{224}}
\put(2937,4062){\ellipse{224}{224}}
\put(5412,4062){\ellipse{224}{224}}
\put(2937,1812){\ellipse{224}{224}}
\put(462,1812){\ellipse{224}{224}}
\thinlines
\path(12,4512)(2487,4512)(2487,2262)
    (12,2262)(12,4512)
\path(2487,4512)(4962,4512)(4962,2262)
    (2487,2262)(2487,4512)
\path(4962,4512)(7212,4512)(7212,2262)
    (4962,2262)(4962,4512)
\path(4962,2262)(7212,2262)(7212,12)
    (4962,12)(4962,2262)
\path(2487,2262)(4962,2262)(4962,12)
    (2487,12)(2487,2262)
\path(12,2262)(2487,2262)(2487,12)
    (12,12)(12,2262)
\thicklines
\path(3387,3500)(3387,3275)
\path(3387,3500)(3387,3275)
\path(3499,3500)(3499,3275)
\path(3499,3500)(3499,3275)
\path(3387,1250)(3387,1025)
\path(3387,1250)(3387,1025)
\path(3499,1250)(3499,1025)
\path(3499,1250)(3499,1025)
\path(1249,3387)(1252,3390)(1259,3396)
    (1270,3406)(1286,3420)(1305,3436)
    (1326,3453)(1347,3469)(1367,3484)
    (1386,3496)(1404,3506)(1420,3513)
    (1434,3518)(1448,3521)(1461,3521)
    (1474,3518)(1487,3514)(1500,3508)
    (1514,3500)(1528,3491)(1542,3480)
    (1558,3468)(1573,3455)(1589,3441)
    (1605,3427)(1622,3413)(1638,3400)
    (1654,3387)(1670,3376)(1686,3366)
    (1702,3357)(1718,3350)(1732,3345)
    (1747,3341)(1763,3338)(1779,3336)
    (1796,3336)(1813,3337)(1831,3339)
    (1850,3342)(1868,3347)(1886,3352)
    (1905,3359)(1923,3366)(1940,3374)
    (1957,3383)(1973,3393)(1989,3403)
    (2004,3414)(2018,3425)(2032,3436)
    (2046,3449)(2061,3462)(2075,3476)
    (2089,3490)(2104,3505)(2118,3520)
    (2132,3536)(2147,3551)(2161,3566)
    (2175,3580)(2189,3593)(2202,3606)
    (2215,3617)(2227,3627)(2239,3636)
    (2251,3643)(2262,3650)(2277,3656)
    (2292,3661)(2307,3663)(2323,3664)
    (2338,3664)(2354,3662)(2369,3658)
    (2384,3653)(2398,3648)(2411,3641)
    (2422,3634)(2433,3627)(2442,3620)
    (2450,3612)(2459,3602)(2467,3591)
    (2474,3580)(2481,3568)(2487,3556)
    (2493,3544)(2500,3532)(2507,3521)
    (2515,3510)(2525,3500)(2534,3491)
    (2545,3482)(2557,3474)(2571,3465)
    (2586,3457)(2602,3449)(2618,3442)
    (2635,3436)(2651,3431)(2666,3427)
    (2680,3425)(2693,3425)(2706,3425)
    (2718,3427)(2730,3431)(2742,3436)
    (2754,3442)(2765,3449)(2777,3457)
    (2787,3465)(2797,3474)(2807,3482)
    (2816,3491)(2824,3500)(2835,3510)
    (2845,3521)(2857,3531)(2868,3541)
    (2881,3551)(2893,3559)(2904,3566)
    (2916,3571)(2926,3574)(2937,3575)
    (2947,3574)(2958,3571)(2969,3566)
    (2981,3559)(2993,3551)(3005,3541)
    (3017,3531)(3028,3521)(3039,3510)
    (3049,3500)(3058,3491)(3067,3482)
    (3076,3473)(3086,3463)(3095,3454)
    (3106,3443)(3116,3433)(3125,3424)
    (3135,3414)(3144,3405)(3153,3396)
    (3162,3387)(3172,3377)(3183,3366)
    (3194,3356)(3206,3346)(3218,3336)
    (3230,3328)(3242,3321)(3253,3316)
    (3264,3313)(3274,3312)(3285,3313)
    (3296,3317)(3308,3323)(3322,3332)
    (3338,3345)(3355,3359)(3370,3372)
    (3381,3382)(3386,3386)(3387,3387)
\path(2712,3387)(2715,3384)(2721,3377)
    (2732,3366)(2746,3351)(2764,3332)
    (2782,3312)(2800,3292)(2817,3273)
    (2832,3255)(2845,3239)(2856,3225)
    (2865,3213)(2872,3201)(2877,3191)
    (2881,3181)(2883,3171)(2885,3162)
    (2885,3153)(2884,3144)(2881,3134)
    (2877,3125)(2872,3115)(2865,3106)
    (2857,3097)(2847,3087)(2836,3078)
    (2824,3068)(2812,3059)(2798,3050)
    (2783,3041)(2768,3031)(2755,3023)
    (2740,3015)(2725,3006)(2709,2997)
    (2692,2988)(2674,2978)(2656,2967)
    (2636,2956)(2617,2944)(2596,2932)
    (2576,2920)(2556,2908)(2536,2895)
    (2516,2882)(2497,2870)(2478,2857)
    (2461,2844)(2444,2832)(2427,2819)
    (2412,2806)(2395,2792)(2378,2777)
    (2362,2761)(2346,2745)(2330,2728)
    (2315,2711)(2300,2694)(2286,2676)
    (2272,2659)(2259,2641)(2247,2624)
    (2235,2607)(2225,2592)(2215,2577)
    (2207,2562)(2199,2549)(2193,2536)
    (2187,2525)(2180,2512)(2175,2500)
    (2171,2488)(2167,2477)(2164,2466)
    (2162,2455)(2162,2444)(2162,2434)
    (2164,2423)(2166,2413)(2170,2403)
    (2175,2394)(2181,2384)(2188,2375)
    (2196,2366)(2205,2356)(2214,2347)
    (2224,2338)(2236,2328)(2248,2318)
    (2261,2307)(2275,2296)(2290,2283)
    (2305,2271)(2321,2257)(2337,2244)
    (2352,2230)(2368,2216)(2383,2202)
    (2398,2188)(2412,2173)(2425,2159)
    (2438,2145)(2449,2131)(2461,2117)
    (2471,2102)(2481,2086)(2491,2070)
    (2501,2053)(2509,2035)(2518,2016)
    (2526,1998)(2533,1978)(2539,1959)
    (2545,1940)(2549,1920)(2553,1901)
    (2556,1883)(2559,1865)(2560,1847)
    (2561,1829)(2562,1812)(2561,1795)
    (2560,1777)(2559,1760)(2557,1741)
    (2554,1723)(2551,1703)(2547,1684)
    (2543,1664)(2538,1643)(2533,1623)
    (2527,1603)(2522,1584)(2516,1564)
    (2510,1546)(2504,1527)(2498,1510)
    (2493,1492)(2487,1475)(2481,1458)
    (2475,1440)(2470,1422)(2464,1404)
    (2458,1386)(2453,1367)(2448,1347)
    (2443,1328)(2438,1309)(2434,1289)
    (2431,1271)(2429,1252)(2427,1234)
    (2426,1217)(2426,1201)(2426,1185)
    (2428,1170)(2430,1156)(2434,1142)
    (2438,1128)(2443,1114)(2450,1100)
    (2457,1086)(2466,1072)(2475,1058)
    (2486,1045)(2497,1032)(2509,1020)
    (2522,1008)(2535,997)(2548,987)
    (2562,978)(2576,970)(2590,962)
    (2604,956)(2618,949)(2634,943)
    (2651,938)(2669,933)(2687,930)
    (2706,926)(2726,924)(2746,923)
    (2765,922)(2785,923)(2804,924)
    (2823,926)(2840,930)(2856,933)
    (2872,938)(2886,943)(2899,949)
    (2914,957)(2927,966)(2940,977)
    (2953,988)(2965,1000)(2977,1013)
    (2988,1027)(2998,1041)(3008,1055)
    (3017,1069)(3026,1082)(3034,1095)
    (3042,1107)(3049,1118)(3057,1129)
    (3064,1140)(3072,1150)(3080,1160)
    (3089,1169)(3098,1178)(3108,1186)
    (3118,1193)(3128,1199)(3138,1204)
    (3149,1208)(3159,1210)(3170,1212)
    (3181,1212)(3192,1212)(3204,1210)
    (3218,1207)(3233,1202)(3251,1196)
    (3270,1189)(3292,1180)(3315,1170)
    (3338,1160)(3358,1151)(3373,1144)
    (3382,1139)(3386,1137)(3387,1137)
\path(3499,1137)(3501,1134)(3505,1129)
    (3513,1119)(3525,1104)(3540,1085)
    (3559,1061)(3581,1035)(3604,1006)
    (3628,977)(3653,947)(3677,919)
    (3700,892)(3722,868)(3743,845)
    (3763,823)(3782,804)(3801,786)
    (3819,769)(3838,754)(3856,739)
    (3874,724)(3895,709)(3917,694)
    (3939,680)(3962,666)(3986,652)
    (4011,639)(4036,626)(4061,614)
    (4087,603)(4113,593)(4138,583)
    (4163,575)(4187,568)(4210,563)
    (4232,559)(4253,555)(4273,554)
    (4291,553)(4308,554)(4324,556)
    (4339,559)(4353,564)(4367,569)
    (4379,577)(4391,585)(4403,595)
    (4413,606)(4423,618)(4432,632)
    (4440,647)(4447,662)(4453,678)
    (4459,695)(4463,712)(4467,730)
    (4470,747)(4472,765)(4473,783)
    (4474,801)(4474,819)(4474,837)
    (4473,856)(4472,875)(4470,895)
    (4468,916)(4464,938)(4460,960)
    (4456,984)(4451,1008)(4445,1033)
    (4439,1058)(4432,1083)(4424,1109)
    (4417,1135)(4408,1160)(4400,1186)
    (4391,1211)(4381,1237)(4372,1262)
    (4362,1287)(4353,1309)(4344,1331)
    (4334,1354)(4324,1378)(4314,1402)
    (4304,1427)(4293,1453)(4282,1479)
    (4271,1506)(4260,1534)(4249,1562)
    (4238,1590)(4228,1618)(4218,1647)
    (4208,1675)(4200,1703)(4191,1730)
    (4184,1757)(4177,1783)(4171,1809)
    (4166,1834)(4162,1858)(4158,1882)
    (4156,1906)(4154,1932)(4153,1958)
    (4153,1983)(4154,2009)(4157,2036)
    (4160,2062)(4164,2088)(4169,2115)
    (4175,2141)(4182,2167)(4190,2193)
    (4199,2218)(4208,2242)(4217,2266)
    (4228,2288)(4238,2309)(4249,2329)
    (4260,2348)(4271,2366)(4283,2382)
    (4294,2398)(4306,2412)(4319,2427)
    (4332,2441)(4346,2455)(4360,2467)
    (4375,2478)(4391,2489)(4406,2499)
    (4423,2507)(4439,2514)(4456,2520)
    (4472,2525)(4488,2528)(4505,2530)
    (4520,2530)(4536,2529)(4551,2527)
    (4565,2524)(4579,2519)(4592,2513)
    (4605,2506)(4618,2498)(4631,2488)
    (4643,2477)(4656,2465)(4669,2451)
    (4682,2436)(4695,2420)(4708,2402)
    (4721,2384)(4734,2365)(4747,2345)
    (4759,2325)(4772,2304)(4784,2284)
    (4795,2264)(4807,2244)(4818,2224)
    (4828,2205)(4839,2187)(4849,2169)
    (4860,2149)(4872,2129)(4884,2109)
    (4896,2089)(4909,2069)(4922,2049)
    (4936,2030)(4950,2010)(4964,1991)
    (4979,1973)(4993,1956)(5007,1940)
    (5022,1925)(5036,1911)(5050,1899)
    (5065,1888)(5079,1877)(5093,1868)
    (5109,1859)(5126,1851)(5145,1844)
    (5165,1839)(5187,1833)(5211,1829)
    (5239,1825)(5268,1822)(5299,1819)
    (5329,1816)(5357,1815)(5380,1813)
    (5397,1813)(5407,1812)(5411,1812)(5412,1812)
\put(980,3000){\makebox(0,0)[lb]{$0$}}
\put(3387,2840){\makebox(0,0)[lb]{$b_1$}}
\put(3387,1362){\makebox(0,0)[lb]{$b_2$}}
\put(5564,1475){\makebox(0,0)[lb]{$y$}}
\put(2599,3612){\makebox(0,0)[lb]{$z$}}
\put(349,3700){\makebox(0,0)[lb]{$x$}}
\end{picture}
}
  \caption{Illustration of the right hand side in \refeq{lemBK}. The bond $b_1$ has been explored first and is found to be $\Z$-vacant and $\T$-vacant. The bond $b_2$ is $r$-equivalent and thus has been explored in the $\Z$ exploration only, it is $\T$-vacant (determined by $b_1$ during the $\T$-exploration), but $\Z$-occupied.}
  \label{figure1}
}
\end{figure}

The exploration process of $\Ctorus(0)$ must be completed at some point,
since the number of bonds within $\Ctorus(0)$ is finite and vertices
turn active only if a bond containing it is explored and is found to be occupied.
We call the result of this exploration process the
{\it $\T$-exploration}.
The cluster $\Ctorus(0)$ consists of all vertices in $\torus$ that are contained in a bond that is $r$-equivalent to a black bond.
However, we have embedded the cluster $\Ctorus(0)$ into
$\Z^d$, which will be useful when we also wish to describe
the related cluster
$\CZ(0)$.

For $\CZ(0)$, the exploration of the cluster is similar, but there are no
gray bonds.
We start with the final configuration of the $\T$-exploration, and
set all vertices that are a common endpoint of a black bond and a gray bond active.
Then we make all gray bonds uncolored again.
From this setting we apply the coloring scheme that colors the uncolored bonds that contain an active vertex.
The coloring scheme is the one for the exploration cluster on $\Zd$, where no gray bonds are created.
Again, we perform this exploration until there are no more active vertices, i.e.,
no more black bonds attached to uncolored bonds. The result is called
the {\it $\Z$-exploration}. In particular, the black bonds in the $\T$-exploration
are a subset of the black bonds in the $\Z$-exploration, which proves that
$\{0\Tconn x\}\subseteq\bigcup_{y\simr x}\{0\Zconn y\}$, and hence $|\Ctorus(0)|\leq |\CZ(0)|$.\footnote{To obtain a coupling for the full percolation configurations on $\torus$ and $\Zd$, we can finally let all bonds that have not been explored be independently occupied with probability $p$, both in $\torus$ and in $\Zd$, independently of each other. However, we do not rely on the coupling of the percolation configuration, but only on the coupling of the clusters $\Ctorus(0)$ and $\CZ(0)$.}

We now show \refeq{lemBK}.
When $\{0\Zconn y\}$ occurs, then picture all $\Z$-occupied paths
from $0$ to $y$ in mind. Since $x\simr y$ and $\{0\Tconn x\}^c$ occurs,
each of these paths $\Z$-connecting $0$ and $y$ should contain a bond which is $\T$-vacant.
Fix such a (self-avoiding) path $\omega\colon0\conn y$ that is $\Z$-occupied, and denote by $b_2$
the first bond that is $\T$-vacant, but $\Z$-occupied,
so that $(0\Zconn \bb_2)\circ
(b_2 \text{ is }\Z-\text{occ.})\circ (\tb_2\Zconn y)$.
Due to our coupling, this implies that there
exists a previously explored bond $b_1$ that is $r$-equivalent to $b_2$,
which is ($\Tbold$- and $\Z$-) vacant. This, in turn, implies that there exists a vertex $z$ that is visited by $\omega$
such that $(z\Zconn \bb_1)$ without using any of the bonds in $\omega$.
Therefore, the event $(0\Zconn z)\circ (z\Zconn \bb_1)\circ (z\Zconn \bb_2)\circ
(b_2 \text{ is }\Z-\text{occ.})\circ (\tb_2\Zconn y)$ occurs (see Figure \ref{figure1}).
\qed
\vskip0.5cm

\section{Previous Results}\label{sec-previous}
In this section we cite results of previous works that will be used in our analysis later on.
We make essential use of results by Borgs, Chayes, van der Hofstad, Slade and Spencer \cite{BCHSS04a,BCHSS04b} for percolation on $\torus$, which we cite in the following two theorems.
\bt[Subcritical phase]\label{thm-SubcriticalPhase}
Under the conditions in Theorem \ref{thm-1}, for $\lambda$ sufficiently small and any $q\ge0$,
    \eq\lbeq{SubcriticalPhase}
    \left(\lambda^{-1}V^{-1/3}+q\right)^{-1}\le\chi_\sT\!\left(p_c(\torus)-\cn^{-1}q\right)\le\left(\lambda^{-1}V^{-1/3}+q/2\right)^{-1}.
    \en
Also, for $p=p_c(\torus)-\cn^{-1}q$ and $\omega\ge1$,
    \eq\lbeq{SubcriticalPhase3}
    \Pbold_{\sT,p}\left(|\Cmax|\ge\frac{\chi_\sT^2(p)}{3600\,\omega}\right)
    \ge\left(1+\frac{36\,\chi_\sT^3(p)}{\omega\,V}\right)^{-1}.
    \en
\et
\noindent
Instead of the upper bound in \refeq{SubcriticalPhase}, we will mainly use the cruder bound
    \eq\lbeq{SubcriticalPhase2}
        \chi_\sT\!\left(p_c(\torus)-\cn^{-1}q\right)\le\frac2q.
    \en

\begin{theorem}[Scaling window]
\label{thm-3}
Assume the conditions in Theorem \ref{thm-1}.
Let $\lambda>0$ and $\Lambda<\infty$.
Then there is a finite positive constant $b$ (depending on $\lambda$ and $\Lambda$) such that, for $p=p_c(\torus)+\Omega^{-1}\epsilon$ with $|\epsilon|\le\Lambda V^{-1/3}$ and $\omega\ge1$,
\eq
\lbeq{CmaxBound}
\Pbold_{\sT, p}\Big(\omega^{-1}V^{2/3}\le|\Cmax|\le\omega V^{2/3}\Big)\ge 1-\frac{b}{\omega}.
\en
\end{theorem}

Theorems \ref{thm-SubcriticalPhase} and \ref{thm-3} have been proven in \cite[Theorems 1.2 and 1.3]{BCHSS04a} subject to the triangle condition.
Using the lace expansion, the triangle condition was established in \cite[Proposition 1.2 and Theorem 1.3]{BCHSS04b} for $d>6$ and sufficient spread-out or $d$ sufficiently large for nearest-neighbor percolation.

Furthermore, we will use two properties of $\Z^d$-percolation in high dimensions formulated in the following two theorems.
\bt[Expected cluster size]\label{prop-gamma}
Under the conditions in Theorem \ref{thm-1},  there exists a positive constant $\Cg$, such that
\eq\lbeq{defCg}
\frac1{\cn\left({p_c(\Z^d)-p}\right)}\le\chi_\sZ(p)\le\frac\Cg{\cn\left({p_c(\Z^d)-p}\right)}\qquad\text{as $p\nearrow p_c(\Z^d)$}.
\en
In other words, the critical exponent $\gamma$ exists and takes on the mean-field value $1$.
\et
\noindent
This theorem is proven in \cite{HS90a} and \cite{AN84}. According to \cite{HS94}, $d\ge19$ is sufficient for the nearest-neighbor model.

We also need (sub)critical bounds on the decay of the connectivity function. These are expressed in the following theorem.
\bt[Bounds on the two-point function]\label{thm-4}
Under the conditions in Theorem \ref{thm-1},  there exist constants $c_{\sss \tau}, C_{\sss \tau}, \cxi, \Cxi>0$ such that,
    \eq
    \lbeq{powerbd}
    \frac{c_{\sss \tau}}{(|x|+1)^{d-2}} \leq \tau_{\sZ,p_c(\Zd)}(x) \leq \frac{C_{\sss \tau}}{(|x|+1)^{d-2}}.
    \en
In other words, the critical exponent $\eta$ exists and takes the value 0.
Furthermore, for any $p<p_c(\Zd)$,
    \eq
    \lbeq{expbd}
    \tau_{\sZ,p}(x) \leq e^{-\frac{\|x\|}{\xi(p)}},
    \en
where the correlation length $\xi(p)$ is defined by
    \eq
    \lbeq{defXi}
    \xi(p)^{-1}=-\lim_{n\to\infty}\frac1n\log\Pbold_{\sZ,p}\big((0,\dots,0)\conn(n,0,\dots,0)\big),
    \en
and satisfies
    \eq
    \lbeq{xibd}
    \cxi\left(p_c(\Z^d)-p\right)^{-1/2}\le\xi(p)\le\Cxi\left(p_c(\Z^d)-p\right)^{-1/2}.
    \en
\et
The power law bound \refeq{powerbd} is due to Hara \cite{Hara05} for the nearest-neighbor case, and to Hara, van der Hofstad and Slade \cite{HHS03} for the spread-out case.
For the exponential bound \refeq{expbd}, see e.g.\ Grimmett \cite[Prop.\ 6.47]{Gri99}.
Hara \cite{Hara90} proves the bound \refeq{xibd}.

\section{The upper bound on the maximal critical cluster}
\label{sec-upperBound}
The following corollary establishes the upper bound on $p_c(\Z^d)$ and the upper bound on
$|\Cmax|$ in Theorem \ref{thm-1}.
For the proof, we first use Proposition \ref{prop-coupling} to obtain that $\chi_\sT(p)\le\chi_\sZ(p)$. Then we use \refeq{defCg} to turn this into relations between $p_c(\torus)$ and $p_c(\Zd)$. Finally, using Theorem \ref{thm-3} we obtain a bound on $|\Cmax|$.
We now present the details of the proof.
\begin{cor}
\label{cor-ubpc}
Under the conditions of Theorem \ref{thm-1} there exists a constant $\Lambda\ge0$ such that, when $r\rightarrow \infty$,
    \eq
    \lbeq{pcZdub}
    p_c(\Z^d)\leq p_c(\torus)+\frac{\Lambda}{\cn} V^{-1/3}.
    \en
Consequently, for $b$ as in Theorem \ref{thm-3}
and all $\omega\ge1$,
    \eq
    \lbeq{Cmaxub}
    \Pbold_{\sT, p_c(\Z^d)}\Big(|\Cmax|\leq \omega V^{2/3}\Big)
    \geq 1-\frac{b}\omega.
    \en
\end{cor}

\vskip0.5cm
\noindent
\proof
By Proposition \ref{prop-coupling},
    \eq\lbeq{ubpc1}
    \chi_{\sT}(p_c(\torus))\leq \chi_{\sZ}(p_c(\torus)).
    \en
When $p_c(\torus) \geq p_c(\Z^d)$, then \refeq{pcZdub}
holds with $\Lambda=0$, so we will next assume that $p_c(\torus) < p_c(\Z^d).$
Using \refeq{defPcT}, \refeq{ubpc1} and \refeq{defCg}, we obtain that
    \eq
    \lambda V^{1/3} \leq \frac{\Cg}{\cn\left(p_c(\Z^d)-p_c(\torus)\right)},
    \en
so that
    \eq
    p_c(\Z^d) \leq p_c(\torus)+\frac{\Cg}{\lambda \cn} V^{-1/3},
    \en
which is \refeq{pcZdub} with $\Lambda =\lambda^{-1}\Cg$.

The bound \refeq{Cmaxub} follows from the fact that, with
$p=p_c(\torus)+\Lambda\cn^{-1}V^{-1/3}\ge p_c(\Zd)$ by \refeq{pcZdub},
    \eq
    \Pbold_{\sT, p_c(\Z^d)}\Big(|\Cmax|\leq \omega V^{2/3}\Big)
    \geq \Pbold_{\sT, p}\Big(|\Cmax|\leq \omega V^{2/3}\Big)
    \geq 1-\frac{b}\omega
    \en
for some constant $b>0$ depending on $\lambda$ and $\Lambda$, and all $\omega\ge1$.
We have used Theorem \ref{thm-3} in the last bound.
\qed
\vskip 0.5cm

\section{The lower bound on the maximal critical cluster}
\label{sec-lowerBound}
In this section, we will bound $\chi_{\sT}(p)-\chi_{\sZ}(p)$ from below.
First we use the results and framework of Section \ref{sec-coupling} to prove such a lower bound in terms of $\widetilde\chi_{\sZ}(p,r)$. Subsequently, assuming \refeq{sharpassumpthm}, we use Theorem \ref{thm-3}, Corollary \ref{cor-ubpc} and the bounds \refeq{SubcriticalPhase2} and \refeq{CmaxBound} to prove Theorem \ref{thm-2}.

Some more work is required if we do not assume \refeq{sharpassumpthm}. We first deduce a bound on $\widetilde\chi_{\sZ}(p,r)$ using the bounds in Theorem \ref{thm-4}. Then we use this bound together with Lemma \ref{lem-BdChiBelow} to show Theorem \ref{thm-1} with the same ingredients as for the proof of Theorem \ref{thm-2}.

\subsection{A lower bound on $\chi_\sT(p)$ in terms of $\chi_\sZ(p)$}
\bl\label{lem-BdChiBelow}
For all $p\in[0,1]$ and $r\ge3$ in the nearest-neighbor model or $r\ge2L+1$ in the spread-out model,
    \eq\lbeq{BdChiBelow}
    \chi_{\sT}(p)
    \geq \chi_{\sZ}(p) \big(1-\chi_{\sZ}(p)\,\widetilde \chi_{\sZ}(p,r)
    -p \,\cn^2\chi_{\sZ}(p)^2\widetilde \chi_{\sZ}(p,r)\big).
    \en
\el
\noindent
\proof
The bound \refeq{BdChiBelow} will be achieved by comparing the two-point functions
on the torus and on $\Z^d$. We will write
$\Pbold=\Pbold_{\sZ,\sT}$
and omit the percolation parameter $p$ from the notation.
Using \refeq{subsetCs}, we write
    \eq
    \tau_{\sT}(x) = \Pbold\Big(\bigcup_{y\in \Z^d: y\simr x} \{0\Zconn y\}\Big)
    -\Pbold\Big(\bigcup_{y\in \Z^d: y\simr x} \{0\Zconn y\}\cap \{0\Tconn x\}^c\Big).
    \en
We further bound, using inclusion-exclusion,
    \eq
    \Pbold\Big(\bigcup_{y\in \Z^d: y\simr x} \{0\Zconn y\}\Big)
    \geq \sum_{y\in \Z^d: y\simr x}\Pbold\big(0\Zconn y\big)
    -\frac 12 \sum_{y_1\neq y_2\in \Z^d: y_1, y_2\simr x}\Pbold\big(0\Zconn y_1, y_2\big),
    \en
so that
    \eq
    \tau_{\sT}(x) \geq  \sum_{y\in \Z^d: y\simr x}\tau_{\sZ}(y)
    -\frac 12 \sum_{y_1\neq y_2\in \Z^d: y_1, y_2\simr x}\shift\Pbold\big(0\Zconn y_1, y_2\big)
    -\Pbold\Big(\bigcup_{y\in \Z^d: y\simr x} \{0\Zconn y\}\cap \{0\Tconn x\}^c\Big).
    \en
Summation over $x\in \torus$ and using that
$\sum_{x\in \torus}\sum_{y\in \Z^d: y\simr x}=\sum_{y\in \Z^d}$ yields that
    \eq
    \lbeq{chi-lb2}
    \chi_{\sT}(p) \geq \chi_{\sZ}(p)
    -\chi_{\sT,1}(p)-\chi_{\sT,2}(p),
    \en
where
    \eqalign
    \chi_{\sT,1}(p)&=\frac 12 \sum_{y_1\neq y_2\in \Z^d: y_1\simr y_2}\Pbold\big(0\Zconn y_1, y_2\big),\\
    \chi_{\sT,2}(p)&=\sum_{x\in\torus}\Pbold\Big(\bigcup_{y\in \Z^d: y\simr x} \{0\Zconn y\}\cap \{0\Tconn x\}^c\Big).
\enalign
Here we use that the sum over $x$ and over $y_1$ and $y_2$ such that $y_1, y_2\simr x$ is the same as the
sum over $y_1$ and $y_2$ such that $y_1\simr y_2$. We are left to bound $\chi_{\sT,1}(p)$ and $\chi_{\sT,2}(p)$.
We start by bounding $\chi_{\sT,1}(p)$. Using the tree-graph inequality (see \cite{AN84}),
    \eq
    \Pbold\big(0\Zconn x, y\big)\le\sum_{z\in\Zd}\tau_\sZ(z)\,\tau_\sZ(x-z)\,\tau_\sZ(y-z),
    \en
(for which the proof easily follows using the BK-inequality \cite{BK85}),
we obtain
    \eqarray
    \nonumber
    \chi_{\sT,1}(p)
    &\leq&
    \frac12\, \sum_{z} \sum_{y_1\neq y_2: y_1\simr y_2}\,\tau_{\sZ}(z)\,\tau_{\sZ}(y_1-z)\,
    \tau_{\sZ}(y_2-z)\\
    &=&
    \frac12\,\chi_\sZ(p)\;\sum_{y_1'\neq y_2':y_1'\simr y'_2}\tau_\sZ(y'_1)\,\tau_\sZ(y'_2).
    \enarray
Here, and in the remainder of the proof, all sums over vertices will be over $\Z^d$ unless written explicitly otherwise.

Since $y_1'\neq y_2'$ and $y_1'\simr y_2'$, we must have that $\|y_1'\|\geq \frac r2$
or $\|y_2'\|\geq \frac r2$. By symmetry, these give the same contributions, so that,
    \eq
    \lbeq{BdChiBelow1}
    \chi_{\sT,1}(p)
    \leq
    \chi_\sZ(p)
    \sum_{y_1'}
    \tau_{\sZ}(y_1')\!
    \sum_{y_2': y_1'\simr y_2', \|y_2'\|
    \geq \frac r2}
    \tau_{\sZ}(y_2')
    \leq \chi_{\sZ}(p)^2 \,\widetilde \chi_{\sZ}(p,r),
    \en
where we recall \refeq{tildechidef}.
We are left to prove that $\chi_{\sT,2}(p)\le p\,\cn^2\chi_\sZ^3(p)\,\widetilde \chi_{\sZ}(p,r)$.

We use \refeq{lemBK}
and note that the right hand side of \refeq{lemBK} does not depend on $x$.
Since
$\sum_{x\in\torus}\sum_{y\in \Z^d: y\simr x}
=\sum_{y\in \Z^d}$, this brings us to
    \eq
    \chi_{\sT,2}(p)\leq \sum_{y\in\Zd} \Pbold\Big(\bigcup_{b_1\neq b_2: b_1\simr b_2}\bigcup_{z\in \Z^d}
    (0\Zconn z)\circ (z\Zconn \bb_1)\circ (z\Zconn \bb_2)\circ
    (b_2 \text{ is }\Z-\text{occ.})\circ (\tb_2\Zconn y)\Big).
    \en
Therefore, by the BK-inequality \cite{BK85},
    \eqalign
    \chi_{\sT,2}(p) &\leq \sum_{y,z}\sum_{b_1\neq b_2: b_1\simr b_2} \Pbold\big((0\Zconn z)\circ (z\Zconn \bb_1)\circ (z\Zconn \bb_2)\circ
    (b_2 \text{ is }\Z-\text{occ.})\circ (\tb_2\Zconn y)\big)\nonumber\\
    &\leq p\sum_{y,z}\sum_{b_1\neq b_2: b_1\simr b_2}\tau_{\sZ}(z)\,
    \tau_{\sZ}(\bb_1-z)\,\tau_{\sZ}(\bb_2-z)\,\tau_{\sZ}(y-\tb_2).
    \enalign
We can perform the sums over $z,y$ to obtain, with $b_i'=b_i-z$,
    \eqalign
    \chi_{\sT,2}(p)
    &\leq p\,\chi_{\sZ}(p)^2
    \sum_{b_1'\neq b_2'\colon b_1'\simr b_2'}
    \tau_{\sZ}(\bb_1')\,\tau_{\sZ}(\bb_2')
    =p\, \cn^2 \chi_{\sZ}(p)^2 \sum_{u\neq v\colon u\simr v} \tau_{\sZ}(u)\,\tau_{\sZ}(v),
    \enalign
where the factor $\cn^2$ arises from the number of choices for $b_1'$ and $b_2'$ for fixed
$\bb_1'$ and $\bb_2'$.
Therefore, by \refeq{BdChiBelow1}, we arrive at the bound
    \eq\lbeq{BdChiBelow2}
    \chi_{\sT,2}(p)\le p\,\cn^2\chi_\sZ^3(p)\,\widetilde \chi_{\sZ}(p,r)
    \en
The bounds \refeq{BdChiBelow1} and \refeq{BdChiBelow2} complete the proof of \refeq{BdChiBelow}.
\qed
\vskip 0.5cm

\subsection{Proof of the main results}
\noindent
\proof[Proof of Theorem \ref{thm-2}]
We assume \refeq{sharpassumpthm} and
take $p=p_c(\Z^d)-{K_1}{\cn^{-1}}V^{-1/3}$ for $K_1$ sufficiently large.
Choose $V$ sufficiently large to ensure that $p>0$.
When $K_1\geq K$, the bound \refeq{sharpassumpthm} still holds.
We obtain from Lemma \ref{lem-BdChiBelow} together with Theorem \ref{prop-gamma} and \refeq{sharpassumpthm} that
    \eq
    \chi_{\sT}(p) \geq K_1^{-1}V^{1/3} \left(1-\Cg C_{\sss K} K_1^{-1}V^{-1/3}-p\,\cn^2\Cg^2C_{\sss K} K_1^{-2}\right)\geq \widetilde c_{\sss K_1} V^{1/3},
    \en
where $\widetilde c_{\sss K_1}$ is chosen appropriately.
Let $K_1$ be so large that $p<p_c(\torus)$, which can be done by \refeq{pcZdub}.
Then, by \refeq{SubcriticalPhase2},
    \eq
    \frac{2}{\cn(p_c(\torus)-p_c(\Zd)+K_1V^{-1/3})}
    \ge\chi_\sT(p)
    \ge \widetilde c_{\sss K_1} V^{1/3},
    \en
so that
    \eq\lbeq{sharpLowerBd2}
    p_c(\Zd)\ge p_c(\torus)+\left(K_1-\frac2{\widetilde c_{\sss K_1}\cn}\right) V^{-1/3},
    \en
which is \refeq{pcZdubd} with $\Lambda=\big(2\,\widetilde c_{\sss K_1}^{-1}-\cn K_1\big)\vee0$.
This, together with \refeq{pcZdub}, permits using Theorem \ref{thm-3}.
By doing so, we obtain that, for $p$ equal to the right hand side of \refeq{sharpLowerBd2},
    \eq
    \Pbold_{\sT,\, p_c(\Z^d)}\Big(\omega^{-1} V^{2/3}\le|\Cmax|\le\omega V^{2/3}\Big)
    \ge\Pbold_{\sT,\, p}\Big(\omega^{-1} V^{2/3}\le|\Cmax|\le\omega V^{2/3}\Big)
    \ge 1-\frac{b}{\omega}.
    \en
This completes the proof of Theorem \ref{thm-2}.
\qed
\vskip0.5cm

\noindent
Unfortunately, we cannot quite prove \refeq{sharpassumpthm}
so we will give cruder upper bounds on $\widetilde \chi_{\sZ}(p,r)$.
This is the content of the following lemma:

\begin{lemma}
\label{lem-tildechibd}
Under the conditions in Theorem \ref{thm-1},
choose $K$ sufficiently large, and let $R= K(\log{V}) (p_c(\Z^d)-p)^{-1/2}$.
Then for all $p\leq p_c(\Z^d)-{K}{\cn^{-1}}V^{-1/3}$,
    \eq
    \widetilde \chi_{\sZ}(p,r)\leq \frac{\Ctc R^2}{V}.
    \en
for some constant $\Ctc>0$.
\end{lemma}
Note that it is here where the power of $\log V$ comes into play.

\proof
For $p\leq p_c(\Z^d)-{K}{\cn^{-1}}V^{-1/3}$,
we bound
    \eq
    \widetilde \chi_{\sZ}(p,r)
    =\sup_{y}\sum_{z\simr y, \|z\|
    \geq \frac r2} \tau_{\sZ,p}(z)
    \leq \sup_{y}\sum_{z\simr y, \|z\|
    \geq R} \tau_{\sZ,p}(z) +\sup_{y}\sum_{z\simr y, \frac r2 \leq  \|z\|
    \leq R} \tau_{\sZ,p}(z).
    \en
We start with the second contribution, for which we
use \refeq{powerbd}. Since $\|z\|\geq \frac r2$, we have that
    \eq
    \tau_{\sZ,p}(z)\le
    \tau_{\sZ,p_c(\Zd)}(z)\le
     \frac{C_{\sss \tau}}{(|z|+1)^{d-2}}\leq \frac{C_1}{V}\sum_{x\in \T_{r,d}}  \frac{C_{\sss \tau}}{(|z+x|+1)^{d-2}}
     \en
for constants $C_{\sss \tau},C_1>0$, where $C_1$ depends on the dimension $d$ only.
Therefore,
    \eq\lbeq{tildechibd1}
    \sup_{y}\shift \sum_{z\simr y: \frac r2 \leq  \|z\|
    \leq R}\shift \tau_{\sZ,p}(z)\leq \frac{C_1}{V} \sup_{y} \sum_{x\in \T_{r,d}}\sum_{ z\simr y: \|z\|\leq R}
    \frac{C_{\sss \tau}}{(|z+x|+1)^{d-2}}
    \leq \frac{C_1}{V} \shift\sum_{z: \|z\|\leq 2R}\frac{C_{\sss \tau}}{(|z|+1)^{d-2}}
    \leq \frac{C_2 R^2}{V},
    \en
where the positive constant $C_2$ depends on $d$ and $L$ only.
For the sum due to $\|z\|\geq R$, we use \refeq{expbd} and \refeq{xibd}
to see that
    \eq\lbeq{tildechibd2}
    \sup_{y}\sum_{z\simr y,\, \|z\|
    \geq R} \tau_{\sZ,p}(z)
    \leq \sup_{y}\sum_{z\simr y,\, \|z\|
    \geq R} \exp\left\{{-\Cxi^{-1}\,\|z\|\,(p_c(\Z^d)-p)^{1/2}}\right\}.
    \en
Since $\|z\|\ge\frac1d(|z_1|+\dots+|z_d|)$ for all $z=(z_1,\dots,z_d)\in\Zd$,
\refeq{tildechibd2} can be further bounded from above by
    \eqalign
    \nn
    &\sum_{z: \|z\|\geq R} \exp\left\{-\Cxi^{-1}\,\|z\|\,(p_c(\Z^d)-p)^{1/2}\right\}
    \\ \nn
    &\qquad\le\left(\sum_{z_1: |z_1|\geq R} \exp\left\{-(d\Cxi)^{-1}\,|z_1|\,\left(p_c(\Z^d)-p\right)^{1/2}\right\}\right)^d\\
    &\qquad\le C_3
    \left(p_c(\Zd)-p\right)^{-d/2}\exp\left\{{-\frac{R}{\Cxi}\left(p_c(\Z^d)-p\right)^{1/2}}\right\},
    \enalign
for some constant $C_3>0$.
Since $R=K(\log{V}) (p_c(\Z^d)-p)^{-1/2}$, the exponential term can be bounded by
    $V^{-K/\Cxi}$.
Furthermore, by our choice of $p$, $\left(p_c(\Zd)-p\right)^{-d/2}\le \left(K\cn^{-1}\right)^{-d/2}\,V^{d/6}$.
Choose $K$ so large that $K/\Cxi-d/6>1$. Then the upper bound is of the order $o(V^{-1})$.
This, together with \refeq{tildechibd1}, proves the claim.
\qed
\vskip0.5cm

\noindent
We next use Lemma \ref{lem-tildechibd} to prove the lower bound in \refeq{pcZdbd}:
\bl\label{lemmaLowerBoundPc}
    Under the conditions in Theorem \ref{thm-1}, there exists a constant $\Lambda\ge0$ such that
    \eq\lbeq{LowerBoundPc}
    p_c(\Z^d)\ge p_c(\torus)-\frac\Lambda\cn V^{-1/3}(\log{V})^{2/3}.
    \en
\el
\proof
By Lemma \ref{lem-tildechibd}, for all $p\le p_c(\Zd)-K \cn^{-1}V^{-1/3}$,
    \eq
    \tilde\chi_\sZ(p,r)\le\frac{\Ctc K^2(\log V)^2}{V\left(p_c(\Zd)-p\right)}.
    \en
With Theorem \ref{prop-gamma}, this can be further bounded as
    \eq\lbeq{LowerBoundPc1}
    \tilde\chi_\sZ(p,r)\le\Ctc K^2\,\frac{(\log V)^2}{V}\,\cn\,\chi_\sZ(p).
    \en
Then, by Lemma \ref{lem-BdChiBelow} and $\chi_\sZ(p)\ge1$,
    \eq\lbeq{LowerBoundPc2}
    \chi_\sT(p)\ge\chi_\sZ(p)\left(1-(1+p\cn^2)\tilde\chi(p,r)\,\chi_\sZ(p)^2\right)
    \en
if $\cn$ and $r$ are sufficiently large.
Combining \refeq{LowerBoundPc1} and \refeq{LowerBoundPc2} yields
    \eq\lbeq{lowerBdChiT}
    \chi_{\sT}(p) \geq \chi_{\sZ}(p) \left(1-(1+p\cn^2)\cn\,\Ctc K^2\,\frac{(\log{V})^2}{V}\,\chi_{\sZ}(p)^3\right).
    \en

Let
    \eq\lbeq{defHatP}
    \hat p:= p_c(\Z^d)-\frac{\hat C}{\cn}V^{-1/3}(\log{V})^{2/3}
    \en
for some (sufficiently large) constant $\hat C>0$.
Depending on $\hat C$, we take $V$ large to ensure that $\hat p>0$.
Then, by Theorem \ref{prop-gamma},
    \eq\lbeq{LowerBoundPc3}
    \chi_\sZ(\hat p)^3\le\frac{\Cg^3}{\cn^3\left(p_c(\Zd)-\hat p\right)^3}=\left(\frac{\Cg}{\hat C}\right)^3V(\log V)^{-2}.
    \en
Substituting \refeq{LowerBoundPc3} into \refeq{lowerBdChiT} for $p=\hat p$, using $\hat p\le1$ and the lower bound in \refeq{defCg} give
    \eq\lbeq{lowerBdChiT2}
    \chi_\sT(\hat p)\ge\left(1-\frac{(1+\cn^2)\cn\,\Ctc K^2\Cg^3}{{\hat C}^3}\right)\frac{1}{\hat C}\,\frac{V^{1/3}}{(\log V)^{2/3}}.
    \en
We make the $\hat C$ in \refeq{defHatP} so large that
    \eq
    \hat c:=\left(1-\frac{(1+\cn^2)\cn\,\Ctc K^2\Cg^3}{{\hat C}^3}\right)\frac{1}{\hat C}>0,
    \en
so that \refeq{lowerBdChiT2} simplifies to
    \eq\lbeq{lowerBdChiT3}
    \chi_\sT(\hat p)\ge\hat c\, {V^{1/3}}{(\log V)^{-2/3}}.
    \en

The quantity
    \eq
    q:=\cn\left(p_c(\torus)-\hat p\right)={\hat C}V^{-1/3}(\log{V})^{2/3}-\cn\left(p_c(\Zd)-p_c(\torus)\right),
    \en
is positive if $V$ is large enough, by \refeq{pcZdub}.
Hence Theorem \ref{thm-SubcriticalPhase} is applicable, and \refeq{SubcriticalPhase2} yields
    \eq\lbeq{lowerBdChiT4}
    \chi_\sT(\hat p)=\chi_\sT\big(p_c(\torus)-q\,\cn^{-1}\big)\le\frac2{\cn\left(p_c(\torus)-\hat p\right)}.
    \en

Merging \refeq{lowerBdChiT3} and \refeq{lowerBdChiT4}, we arrive at
    \eq
    p_c(\torus)-p_c(\Z^d)\leq \frac1\cn\left[\frac2{\hat c}-{\hat C}\right] V^{-1/3}(\log{V})^{2/3},
    \en
which is \refeq{LowerBoundPc} with $\Lambda=(2\hat c^{-1}-{\hat C})\vee0$.
\qed
\vskip0.5cm

\noindent
\bc\label{cor-lbCmax}
Under the conditions in Theorem \ref{thm-1}, there exists a constant $C>0$ such that, for all $\omega_1\ge C$,
    \eq
    \Pbold_{\sT,p_c(\Zd)}\left(|\Cmax| \geq \omega_1^{-1} (\log{V})^{-4/3}\,V^{2/3}\right)\ge\left(1+\frac{ 120^{3/2}\cdot288}{\omega_1^{3/2}(\log V)^{2}}\right)^{-1}.
    \en
\ec
\proof
Take $V$ so large that $\lambda^{-1} \le \sqrt{\omega_1}\,120^{-1} (\log V)^{2/3}$, and let
    \eq
    \hat p = p_c(\torus) - \sqrt{\omega_1}\,120^{-1} \cn^{-1} V^{-1/3} (\log V)^{2/3}.
    \en
Then, by \refeq{SubcriticalPhase},
\eq
    \chi_{\sT}^2(\hat p)
    \ge\left(\frac1{\lambda^{-1}V^{-1/3}+\sqrt{\omega_1}\,120^{-1}V^{-1/3}(\log V)^{2/3}}\right)^2
    \ge3600\,\omega_1^{-1}(\log{V})^{-4/3}\,V^{2/3}.
\en
This enables the bound
\eq\lbeq{lowerBdCmax}
    \Pbold_{\sT,p_c(\Zd)}\left(|\Cmax| \geq \omega_1^{-1}(\log{V})^{-4/3}\,V^{2/3}\right)
    \ge\Pbold_{\sT,p_c(\Zd)}\left(|\Cmax| \geq \frac{\chi_\sT^2\left(\hat p\right)}{3600}\right).
\en
By Lemma \ref{lemmaLowerBoundPc}, $\hat p\le p_c(\Zd)$ for $\omega_1\ge C$, and $C>0$ large enough.
Thus, we use \refeq{SubcriticalPhase3} and \refeq{SubcriticalPhase2} to bound \refeq{lowerBdCmax} further from below by
\eq
   \Pbold_{\sT,\hat p}\left(|\Cmax| \geq \frac{\chi_\sT^2\left(\hat p\right)}{3600}\right)
   \ge \left(1+\frac{36\,\chi^3_\sT\left(\hat p\right)}{V}\right)^{-1}
   \ge \left(1+\frac{36\cdot2^3\cdot120^{3}}{\omega_1^{3/2}(\log V)^{2}}\right)^{-1}.
\en
\qed
\vskip0.5cm

Combining our results from Sections \ref{sec-upperBound} and \ref{sec-lowerBound}, we finally prove Theorem \ref{thm-1}:
\proof[Proof of Theorem \ref{thm-1}]
By Corollaries \ref{cor-ubpc} and \ref{cor-lbCmax}, for $\omega_1\ge C$ for some sufficiently large $C$ and $\omega_2\ge1$,
    \eq
    \Pbold_{\sT, p_c(\Z^d)}\Big(\omega_1^{-1}(\log{V})^{-4/3}\,V^{2/3}
    \leq |\Cmax|\leq \omega V^{2/3}\Big)
    \geq 1-\left(\frac{\frac{120^{3}\cdot288}{\omega_1^{3/2}(\log V)^{2}}}{1+\frac{120^{3}\cdot288}{\omega_1^{3/2}(\log V)^{2}}}\right)-\frac{b_2}\omega_2,
    \en
where the term in brackets on the right hand side vanishes for $V\to\infty$.
Then $b_1$ in \refeq{Cmaxbd} can be taken as $120^{3}\cdot288$.
This proves Theorem \ref{thm-1}.
\qed
\vskip0.5cm

\subsection{Discussion of \refeq{sharpassumpthm}}
At the end of Section \ref{subsect-results} we argued why we believe that \refeq{sharpassumpthm} holds.
Another approach to \refeq{sharpassumpthm} is to split the sum over $z$ in \refeq{tildechidef}.
Note that, for $p=p_c(\Zd)-K_1\cn^{-1}V^{-1/3}$, the sum due to $\|z\|\leq K_1 V^{1/6}$ can be bounded,
for any $K_1>0$, as
    \eq
    \sup_{y}\sum_{z\simr y, \frac r2\leq \|z\|
    \leq K_1 V^{1/6}} \tau_{\sZ}(z)
    \leq C\sup_{y}\sum_{z\simr y, \frac r2\leq \|z\|
    \leq K_1 V^{1/6}} (\|z\|+1)^{-(d-2)}
    \leq C  K_1^2 V^{-2/3}.
    \en
Therefore, we are left to give a bound on the contribution from $\|z\|\geq K_1 V^{1/6}$.
The restriction $\|z\|\geq K_1 V^{1/6}$ is equivalent to $\|z\|\geq C_{\sss K,K_1} \xi(p),$ where $\xi(p)$ denotes the correlation length. Indeed, $\xi(p)$ is comparable in size to $(p_c(\Z^d)-p)^{1/2}$ as proven by Hara \cite{Hara90} (see also \refeq{xibd}), and the constant $C_{\sss K,K_1}$ can be made arbitrarily large by taking $K_1$ large.
This contribution could be bounded by investigating $\tau_{\sZ}(z)$ for $\|z\|\geq C_{\sss K,K_1} \xi(p).$

\section{The role of boundary conditions}\label{sec-boundaryConditions}
In this section, we discuss the impact of boundary conditions on the geometry of the largest critical cluster.
For the $d$-dimensional box $\{-\lfloor r/2\rfloor,\dots,\lceil r/2\rceil-1\}^d$, we write $\Br$ if we consider it with free boundary conditions,
and we write $\torus$ if it is equipped with periodic boundary conditions.
We fix $p=p_c(\Zd)$ and further omit this subscript.
Furthermore, we write $C$ for a positive constant, whose value may change from line to line.

With Michael Aizenman \cite{Aize05}, we have discussed the role of
boundary conditions for critical percolation above the upper critical dimension.
We will summarize the consequences of Theorems \ref{thm-1} and \ref{thm-2} in this discussion now.
Assume that the conditions in Theorem \ref{thm-1} are satisfied.
Let $X_1, X_2, X_3$ and $X_4$ be 4 uniformly chosen vertices in $\Br$. Then, Aizenman
notices that, with bulk boundary conditions,
    \eq
    \Pbold_{\sZ}(X_1\conn X_3\mid X_1\conn X_2, X_3\conn X_4)
    \rightarrow 0,
    \en
as the width of the torus $r$ tends to infinity. Indeed,
    \eq
    \Pbold_{\sZ}(X_1\conn X_3\mid X_1\conn X_2, X_3\conn X_4)
    =\frac{\Pbold_{\sZ}(X_1\conn X_2,X_3,X_4)}{\Pbold_{\sZ}(X_1\conn X_2,X_3\conn X_4)}.
    \en
We can compute
    \eq
    \Pbold_{\sZ}(X_1\conn X_2,X_3,X_4)=V^{-4} \sum_{x\in \Br} {\mathbb E}_\sZ[|C_{\sZ}(x,r)|^3],
    \en
where $C_{\sZ}(x,r)$ is the set of vertices $y\in\Br$ for which $x\Zconn y$,
and the right hand side will be bounded from above by the following lemma.
\bl\label{lem-treeGraph}
Under the conditions of Theorem \ref{thm-1}, for $p=p_c(\Zd)$ and $r\ge3\vee (2L+1)$,
    \eq
    {\mathbb E}_\sZ[|C_{\sZ}(x,r)|^3]\leq Cr^{10},
    \qquad\text{for all $x\in\Br$.}
    \en
\el
\noindent
The proof will be given at the end of this section.

On the other hand,
    \eqalign
    \Pbold_{\sZ}(X_1\conn X_2, X_3\conn X_4)
    &=V^{-4}\sum_{x,y,u,v\in \Br} {\mathbb P}_\sZ(x\Zconn u,\, y\Zconn v).\nn
    \enalign
By the FKG-inequality, for all $x,y,u,v\in\Br$,
    \eq
    {\mathbb P}_\sZ(x\Zconn u, y\Zconn v)
    \geq  {\mathbb P}_\sZ(x\Zconn u)\,{\mathbb P}_\sZ(y\Zconn v),
    \en
so that
    \eqalign
    \Pbold_{\sZ}(X_1\conn X_2, X_3\conn X_4)
    &\geq V^{-4}\left(\sum_{x,u\in \Br}{\mathbb P}_\sZ(x\Zconn u)\right)^2.
    \enalign
For fixed $x$, by \refeq{powerbd},
    \eq
    \sum_{u\in \Br}{\mathbb P}_\sZ(x\Zconn u)
    \ge \sum_{u\in \Br}\frac{c_{\tau}}{(|x-u|+1)^{d-2}}
    \ge Cr^2.
    \en
Summing over $x$ gives an extra factor $V$.
We obtain that
    \eqalign
    \Pbold_{\sZ}(X_1\conn X_2, X_3\conn X_4)
    &\geq C V^{-2}r^{4}.
    \enalign
Therefore, when $d>6$,
    \eqalign\nn
    \Pbold_{\sZ}(X_1\conn X_3\mid X_1\conn X_2, X_3\conn X_4)
    &=\frac{\Pbold_{\sZ}(X_1\conn X_2,X_3,X_4)}{\Pbold_{\sZ}(X_1\conn X_2,X_3\conn X_4)}\\
    \lbeq{bc1}
    &\leq\frac{V^{-3} r^{10}}{C V^{-2} r^{4}}= \frac{r^{6}}{CV}\rightarrow 0.
    \enalign
All this changes when we consider the torus with periodic boundary conditions and we assume that
    \eq\lbeq{bc5}
    \chi_{\sT}(p_c(\Z^d))=\Theta\big(V^{1/3}\big).
    \en
Note that \refeq{bc5} is a consequence of \refeq{pcZdubd}, which follows from \refeq{sharpassumpthm}, and Theorem \ref{thm-SubcriticalPhase}.
In this case,
    \eq
    \Pbold_{\sT}(X_1\conn X_2,X_3,X_4)\geq \Pbold_{\sT}(X_1, X_2,X_3,X_4 \in \Cmax).
    \en
Thus, for $\omega\ge1$ sufficiently large,
    \eqalign\lbeq{bc2}
    \Pbold_{\sT}(X_1\conn X_2,X_3,X_4)&\geq \Pbold_{\sT}\Big(X_1, X_2,X_3,X_4 \in \Cmax,
    |\Cmax|\geq \frac 1\omega V^{2/3}\Big)\nn\\
    &\geq \omega^{-4} V^{-4/3}
    \Pbold_{\sT}\Big(|\Cmax|\geq \frac 1\omega V^{2/3}\Big)
    \geq \frac 12\omega^{-4} V^{-4/3}.
    \enalign
On the other hand, we have that
    \eq
    \Pbold_{\sT}(X_1\conn X_2, X_3\conn X_4)
    \leq \Pbold_{\sT}\big((X_1\conn X_2)\circ (X_3\conn X_4)\big)+\Pbold_{\sT}(X_1\conn X_2,X_3,X_4),
    \en
and, by the BK-inequality,
    \eq
    \Pbold_{\sT}\big((X_1\conn X_2)\circ (X_3\conn X_4)\big)
    \leq \Pbold_{\sT}(X_1\conn X_2)\,\Pbold_{\sT}(X_3\conn X_4)
    =V^{-2} \chi_{\sT}(p_c(\Z^d))^2\leq C V^{-4/3}.
    \en
Therefore, assuming \refeq{bc5}, we obtain
    \eq\lbeq{bc6}
    \limsup_{V\rightarrow \infty} \Pbold_{\sT}(X_1\conn X_3\mid X_1\conn X_2, X_3\conn X_4)
    >0.
    \en
The difference between \refeq{bc1} and \refeq{bc6} was conjectured by Aizenman \cite{Aize05}.
The obvious conclusion is that boundary conditions play a crucial role
for high-dimensional percolation on finite cubes.

We do not know that \refeq{bc5} holds, so now we will investigate the changes using the results in Theorem \ref{thm-1} in the above discussion. Thus, we will use that, with high probability,
    \eq\lbeq{bc3}
    V^{2/3}(\log{V})^{-4/3}\leq
    |\Cmax| \leq \omega V^{2/3},
    \en
and
    \eq
    \frac{1}{\omega}V^{1/3}(\log{V})^{-2/3}\leq \chi_{\sT}(p_c(\Z^d))\leq \omega V^{1/3}.
    \en
We will see that the conclusion weakens.
Indeed, by \refeq{bc1},
    \eq
    \Pbold_{\sZ}(X_1\conn X_3\mid X_1\conn X_2, X_3\conn X_4)
    \leq{Cr^{6}}{V^{-1}}=C r^{6-d}\rightarrow 0,
    \en
where the convergence is as an inverse power of $r$, while by an argument as in \refeq{bc2}, now using \refeq{bc3},
    \eq
    \Pbold_{\sT}(X_1\conn X_3\mid X_1\conn X_2, X_3\conn X_4)
    \geq C (\log{V})^{-\frac{16}3},
    \en
which only converges to zero as a power of $\log{r}$. Therefore, the main conclusion
that boundary conditions play an essential role is preserved.

\vspace{.5cm}
We have argued that the largest critical cluster with bulk boundary conditions is much smaller than the one with periodic boundary conditions. We will now argue that,
under the condition of Theorem \ref{thm-1}, critical percolation clusters on the periodic torus $\torus$ are similar to percolation clusters on a finite box with bulk boundary conditions, where the box has width $V^{1/6}=r^{d/6}\gg r$.
Here we rely on the coupling in Proposition \ref{prop-coupling}.
In particular, when the origin $0$ is $\Tbold$-connected to a uniformly chosen point $X$, then, with high probability, there is no $\Z$-connection at distance $o(V^{1/6})$ from $0$ to a point that is $r$-equivalent to $X$.

This will be illustrated by the following calculation.
Assume for the moment that $\chi_\sT(p_c(\Zd))=\Theta(V^{1/3})$. This follows from our assumption \refeq{sharpassumpthm}. Choose the vertex $X$ uniformly from the torus $\torus$. Then, for any $\vep>0$,
    \eqalign
    \nonumber
    &\Pbold_{\sZ,\sT}\big(\exists y\in\Zd\colon y\simr X, |y|\le\vep V^{1/6}, 0\Zconn y \mid 0\Tconn X\big)\\
    \lbeq{bc4}
    &\quad\le \frac{\Pbold_{\sZ,\sT}\big(\exists y\in\Zd\colon y\simr X, |y|\le\vep V^{1/6}, 0\Zconn y\big)}{\Pbold_{\sZ,\sT}\big(0\Tconn X\big)}.
    \enalign
For the denominator, we rewrite
    \eq
    \Pbold_{\sZ,\sT}\big(0\Tconn X\big)
    = V^{-1}\chi_\sT(p_c(\Zd))
    \ge CV^{-2/3},
    \en
whereas the numerator in \refeq{bc4} is bounded from above by
    \eq
    \sum_{y\in\Zd}\Pbold_{\sZ,\sT}\!\left(y\simr X, |y|\le\vep V^{1/6}, 0\Zconn y\right)
    \le \frac1V\sum_{y: |y|\le\vep V^{1/6}}\frac1{\left(|y|+1\right)^{d-2}}
    \le C \vep^2V^{-2/3}.
    \en
Thus,
    \eq
    \Pbold_{\sZ,\sT}\big(\exists y\in\Zd\colon y\simr X, |y|\le\vep V^{1/6}, 0\Zconn y \mid 0\Tconn X\big)
    \le C\vep^2.
    \en

We have seen that, under bulk boundary conditions, $|\Cmax|$ is of
the order $\Theta(r^4)$, whereas under periodic boundary conditions,
it is of the order $\Theta(V^{2/3})$. Thus, the maximal critical
percolation cluster on a high dimensional torus is $\Theta(R^4)$
with $R=V^{1/6}\gg r$, so that $\Cmax$ is 4-dimensional, but now in
a box of width $R$. This suggests that percolation on a box
$B_{r,d}$ with periodic boundary conditions is similar to
percolation on the larger box $B_{R,d}$ under bulk boundary
conditions, with $R=V^{1/6}\gg r$.

\vspace{.5cm}
Without assuming \refeq{sharpassumpthm}, the lower bound on the denominator is only $CV^{-2/3}(\log V)^{2/3}$, thus we obtain the weaker bound
    \eq
    \Pbold_{\sZ}\!\left(\exists y\in\Zd\colon y\simr X, |y|\le\vep V^{1/6}(\log V)^{-1/3}, 0\Zconn y \mid 0\Tconn X\right)\le\vep^2.
    \en
The conclusion that occupied paths are long is preserved.
\vspace{.5cm}

We conclude this section with the proof of Lemma \ref{lem-treeGraph}.
\proof[Proof of Lemma \ref{lem-treeGraph}]
For all $x\in\Br$,
    \eqalign
    \lbeq{treeGraph1}
    &{\mathbb E}_\sZ[|C_{\sZ}(x,r)|^3]
    =\sum_{s,t,u\in\Br}\Pbold_\sZ\left(x\Zconn s\Zconn t\Zconn u\right)\\
    \nonumber
    &\qquad\le3\sum_{\substack{s,t,u\in\Br\\v,w\in\Zd}}\Pbold_\sZ\left((x\Zconn v)\circ (v\Zconn s)\circ (v\Zconn w)\circ (w\Zconn t)\circ (w\Zconn u)\right).
    \enalign
Using the BK-inequality, this can be further bounded from above by
    \eq\lbeq{treeGraph2}
    3\bigg(\sup_{v\in\Zd}\sum_{s\in B_{r,d}}\tau_\sZ(s-v)\bigg)^2\sum_{u\in B_{r,d}}\tau_\sZ^{\ast3}(x-u),
    \en
where $\tau^{\ast3}_\sZ$ denotes the threefold convolution of $\tau_\sZ$.
We begin to bound the expression in parenthesis. Fix $v\in\Zd$.
If the distance between $v$ and the box $B_{r,d}$ is larger than $r$, then $\tau_\sZ(s-v)\le Cr^{-(d-2)}$ for all $s\in B_{r,d}$, by \refeq{powerbd}. Hence, in this case,
    \eq
    \sum_{s\in B_{r,d}}\tau_\sZ(s-v)\le Cr^2.
    \en
Otherwise, $\|s-v\|\le2r$ for all $s\in B_{r,d}$, and therefore, by \refeq{tildechibd1},
    \eq
    \sum_{s\in B_{r,d}}\tau_\sZ(s-v)
    \le \sum_{z\in B_{2r,d}}\tau_\sZ(z)
    \le Cr^2.
    \en
Using \cite[Prop.\ 1.7 (i)]{HHS03} for $d>6$,
we see that, for all $z\in\Zd$, the upper bound in \refeq{powerbd}
implies that
    \eq
    \tau_\sZ^{\ast2}(z)\le\frac{C}{(|z|+1)^{d-4}},
    \en
which in turn implies, when $d>6$, so that $(d-2)+(d-4)>d$,
    \eq
    \tau_\sZ^{\ast3}(z)=(\tau_\sZ\ast\tau_\sZ^{\ast2})(z)\le\frac{C}{(|z|+1)^{d-6}}.
    \en
Thus we obtain, for $x\in B_{r,d}$, 
    \eq\lbeq{treeGraph3}
    \sum_{u\in B_{r,d}}\tau_\sZ^{\ast3}(x-u)
    \le \sum_{z\in B_{2r,d}}\tau_\sZ^{\ast3}(z)
    \le \sum_{z\in B_{2r,d}}\frac C{(|z|+1)^{d-6}}
    \le Cr^6.
    \en
The combination of the bounds \refeq{treeGraph1}--\refeq{treeGraph3} yields the desired upper bound $Cr^{10}$.
\qed

\vspace{.5cm}
\noindent
\textbf{Acknowledgment.}
This work was supported in part by the Netherlands Organisation for Scientific Research
(NWO).
The research of RvdH was performed in part while visiting Microsoft Research in the summer 2004. We thank Christian Borgs, Jennifer Chayes and Gordon Slade for valuable discussions at the start of this project, and Michael Aizenman for useful discussions on the relations between scaling limits and IIC's, as well as the role of boundary conditions for high-dimensional percolation as described in Section \ref{sec-boundaryConditions}.
We thank Jeffrey Steif for pointing our attention to \cite{BenjaSchra96},
and the referee for numerous suggestions to improve the presentation.

\end{document}